\documentclass[12pt]{article}
\usepackage{calrsfs}
\usepackage{amsmath}
\usepackage{amsfonts}
\usepackage{latexsym}
\usepackage{amssymb}
\usepackage{amsthm}
\usepackage[normalem]{ulem}
\usepackage{csquotes}
\usepackage{lineno,xcolor}
\usepackage{url}
\usepackage[normalem]{ulem}
\usepackage{filecontents}
\usepackage{scalefnt}
\usepackage{tikz}
\usepackage{fancyhdr}
\usetikzlibrary{arrows}
\usetikzlibrary{positioning}

\usepackage{lineno,xcolor}
\definecolor{amaranth}{rgb}{0.9, 0.17, 0.31}
\definecolor{bluegray}{rgb}{0.4, 0.6, 0.8}

\makeatletter

\usepackage{longtable}

\newtheorem*{maintheorem*}{Main Theorem}
\newtheorem{theorem}{Theorem}[section]
\newtheorem{proposition}[theorem]{Proposition}
\newtheorem{corollary}[theorem]{Corollary}
\newtheorem{lemma}[theorem]{Lemma}

\newtheorem*{theorem*}{Theorem}
\newtheorem{remark}[theorem]{Remark}

\newtheorem*{example*}{Example}
\newtheorem*{conjecture*}{Conjecture}
\def\1{\mathbf 1}

\def\0{\mathbf 0}

\def\cC{\mathcal C}

\def\cL{\mathcal L}

\def\cP{\mathcal P}
\def\cX{\mathcal X}
\def\cV{\mathcal V}

\def\cR{\mathcal R}
\def\cS{\mathcal S}

\def\cV{\mathcal V}
\def\cY{\mathcal Y}

\def\mI{{\,\mathrm I}\,}

\def\<{\langle}
\def\>{\rangle}
\newcommand\comment[1]{}

\newcommand*{\shifttext}[2]{
	\settowidth{\@tempdima}{#2}
	\makebox[\@tempdima]{\hspace*{#1}#2}
}

\newcommand\redsout{\bgroup\markoverwith{\textcolor{amaranth}{\rule[0.5ex]{2pt}{0.4pt}}}\ULon}

\newcommand{\bluemunsell}[1]{\textcolor[rgb]{0.0, 0.5, 0.69}{#1}}

  \makeatletter
\def\@fnsymbol#1{\ensuremath{\ifcase#1\or *\or \dagger\or \ddagger\or
   \mathsection\or \mathparagraph\or \|\or **\or \dagger\dagger
   \or \ddagger\ddagger \else\@ctrerr\fi}}
    \makeatother

\title{The association scheme on the set of flags \\ of a finite generalized quadrangle}

\author{Francesco Colangelo\thanks{Francesco Colangelo: francesco.colangelo@unibas.it\hfill\newline \hspace*{1.4em}
		Alessandro Siciliano: alessandro.siciliano@unibas.it\hfill\newline \hspace*{1.4em}
		Universit\`{a} degli Studi della Basilicata, DiMIE, Viale dell'Ateneo Lucano 10, 85100 Potenza, Italy}	\and Giusy Monzillo \footnote{Giusy Monzillo: giusy.monzillo@famnit.upr.si\hfill\newline \hspace*{1.4em}
		University of Primorska, UP FAMNIT, Glagolj\v aska 8, 6000 Koper, Slovenia}
	\and
	Alessandro Siciliano${}^{\dagger}$
}
\date{}

\begin{document}
	

	\maketitle
	
	\thispagestyle{fancy}
	\fancyhf{}
	\renewcommand{\headrulewidth}{0pt}
	
	\begin{abstract} 
		In this paper, the association scheme defined on the flags of a finite generalized quadrangle is considered. All possible fusions of this scheme are listed, and a full description for those of classes 2 and 3 is given. 
		
	Furthermore, it is showed that an association scheme with appropriate parameters must arise from the flags of a generalized quadrangle. The same is done for one of its 4-class symmetric fusion.
	\end{abstract}

	\section{Introduction}	
	
	An \emph{association scheme} is a pair $\cX=(X,\cR)$ where $X$ is a finite set  and $\cR=\{R_i\}_{i\in I}$ is a collection of binary relations on $X$ satisfying the following properties:
	\begin{enumerate}
		\item[(AS1)] $\cR$ is a partition of $X \times X$.
		\item[(AS2)] The diagonal relation $R_0=\{(x,x):x \in X\} $ is in $\cR$.
		\item[(AS3)] For each $i\in I$, there exists an $i^*\in I$ such that $R_{i^*}=\{(y,x)\colon (x,y)\in R_i\}$  is in $\cR$.
		\item[(AS4)] For each $i,j,k\in I$, there exist constants $p_{ij}^k$, such that if $(x,y)\in R_k$, then there are $p_{ij}^k$ vertices $z$ such that $(x,z)\in R_i$ and $(z,y)\in R_j$. The $p_{ij}^k$ are called {\em intersection numbers}.
	\end{enumerate}

	The cardinality of $X$ is called the \emph{order} of the scheme $\cX$ and that of $I^*=I\setminus \{0\}$ is called the \emph{class} of $\cX$. The relations $R_i \in \cR$ are called \emph{basis relations}, and the digraphs $(X,R_{i}$) \emph{basis digraphs} of $\cX$. It follows from axiom (AS4) that every basis digraph  $(X,R_{i}$) is regular. The {\em valency} of $R_{i}$ is the outdegree of $(X,R_{i})$ and will be denoted by $\eta_{i}$. For every $x\in X$, we set $R_i(x)=\{y \in \Omega	: (x,y) \in R_i  \}.$ An association scheme is said to be \emph{symmetric} if each relation $R_i$ is equal to its opposite $R_{i^*}$;  it is \emph{commutative} if $p^k_{ij}=p^k_{ji}$, for all $i,j,k \in I$.
	
	An association scheme  $\cX=(X,\cR)$ is said to be \emph{thin} if $\eta_{i}=1$ for all $i \in I^*$. For any two given basis relations in $\cR$, the set 
	\[
	R_i R_j = \{R_k : p^k_{ij} \neq 0\}
	\]
	is called the \emph{complex product} of $R_i$ and $R_j$.
	
	It is known that for any thin scheme  $\cX=(X,\cR)$, the set $\cR$ endowed with the complex product is a multiplicative group, whose identity element is $R_0$.

	A \emph{$\{0,i\}$-clique} of $\cX$, with  $i\neq 0$ and $R_i\in\cR$ a symmetric relation, is any clique in the graph  $(X,R_i)$, i.e. any complete subgraph of $(X,R_i)$; a $\{0,i\}$-clique is said to be \emph{maximal} if it is not contained in a larger $\{0,i\}$-clique. 
	
	A union of basis relations of $\cX=(X,\cR)$ which is an equivalence relation on $X$ is called a \emph{parabolic}  of $\cX$. The set of the equivalence classes of a parabolic $e$ is denoted by $X/e$. The parabolics $R_0$ and $X\times X$ are called {\em trivial parabolics}.   A scheme $\cX$ is said to be {\em primitive} if the only parabolics are the trivial ones, and it is called \emph{imprimitive} otherwise.
	
	Let  $\cX=(X,\cR)$ be an imprimitive association scheme with a non-trivial parabolic $e$. By following  \cite[Section 3.1.2]{cp}, it is possible to  construct a  scheme on $X/e$, which is called the {\em quotient scheme} of $\cX$ {\em modulo the parabolic} $e$; this scheme is denoted by $\cX_{X/e}$.
	
	Let $\cX=(X,\{R_i\}_{i\in I}) $ and $\cX'=(X',\{R'_i\}_{i\in I}) $ be association schemes. A bijection 
	\begin{center}
		$\phi: R_i \in \cR \rightarrow R'_{i'} \in \cR' $
	\end{center}
	is called an \emph{algebraic isomorphism} from $\cX$ to $\cX'$ if 
	\begin{center}
		$p^k_{ij}={p}^{k'}_{i'j'}$ \ \ \ for all $i,j,k \in I$. 
	\end{center}
	If such bijection  exists then $\cX$ and $\cX'$ are said to be \emph{algebraically isomorphic}. If $\cX$ and $\cX'$ are algebraically isomorphic then every algebraic isomorphism induces a bijection between the set of parabolics of $\cX$ and the set of parabolics of $\cX'$ \cite[Prop. 2.3.25]{cp}.
	
	The reader is referred to \cite{bi,cam,cp} for additional information on association schemes.

	A (finite) {\em generalized quadrangle} (GQ) is an incidence structure $\cS=(\cP,\cL, {\rm I})$ where  $\cP$ and $\cL$ are disjoint non-empty sets of objects called {\em points} and {\em lines}, respectively, and $\rm I$ is a symmetric point-line incidence relation  satisfying the following axioms:
	\begin{itemize} 
		\item[ (GQ1)] Each point is incident with $t+1$ lines, and two distinct points are incident with at most one line.
		\item[(GQ2)]  Each line is incident with $s+1$ points, and two distinct lines are incident with at most one point.
		\item[(GQ3)] If $p$ is a point and $L$ is a line not incident with $p$, then there is a unique pair $(q, M)\in \cP\times \cL$ such that $p\, {\rm I}\, M\, {\rm I}\, q\, {\rm I}\,  L$.
	\end{itemize}
	
	The integers $s$ and $t$ are the {\em parameters} of the $\mathrm{GQ}$, and $\cS$ is said to have {\em order} $(s,t)$.  If $\cS$ has order $(s,t)$, then $|\cP| = (s + 1)(st + 1)$ and $|\cL| = (t + 1)(st + 1)$. 
	
	Two distinct points $p$ and $q$ are said to be \emph{collinear} on the line $L$ if there is a (unique) line $L$ incident with both $p$ and $q$; for distinct lines $L$ and $M$, we say that $L$ \emph{intersects} $M$ at the point $r$, and we write $L\cap M=\{r\}$ if there is a (unique) point $r$ incident with both $L$ and $M$.
	
	In any $\mathrm{GQ}(s,t)$ 	there is the so-called {\em point-line duality}: in any definition or theorem the words \enquote{point} and \enquote{line}, \enquote{collinear} and \enquote{intersecting}, as well as the parameters, are interchanged. Therefore,  the incidence structure $\cS^D=(\cL,\cP,{\rm I}^D)$, with ${\rm I}^D={\rm I}$, is a generalized quadrangle of order $(t,s)$, called the {\em dual} of $\cS$. 
	
	For more details on generalized quadrangles, the reader is referred to \cite{pt}.
	
	It is known that the points of a generalized quadrangle under the relation of collinearity form a strongly regular graph, called the {\em point-graph} of the quadrangle, which actually gives rise to a 2-class association scheme.
	
	By using the geometry of generalized quadrangles which satisfy prescribed properties, it is possible to construct association schemes with more than two classes. 
	Payne in \cite{payne} constructed a $3$-class association scheme $\cX$ starting from a generalized quadrangle with a quasi-regular point. Subsequently, Hobart and Payne in \cite{HobPay} proved that an association scheme having the same parameters as $\cX$ and satisfying an assumption about certain maximal cliques is necessarily the scheme $\cX$.
	In \cite{gl}, Ghinelli and L\"owe define a $4$-class association scheme on the points of a generalized quadrangle with a regular point, and they characterize the scheme by its parameters.
	Penttila and Williford \cite{pw} constructed an infinite family of 4-class association schemes starting from a generalized quadrangle with a doubly subtended subquadrangle. These schemes have been characterized by their parameters in \cite{ms}. Similarly, it has been done in \cite{mon} for a $4$-class scheme constructed by van Dam, Martin and Muzychuk from a generalized quadrangle with a hemisystem \cite{vdmm}.
	
	\bigskip 
	In the spirit of \cite{gl,HobPay,mon,ms,payne} we consider the scheme on the set of the flags (incident point-line  pairs) of a generalized quadrangle. In Section \ref{sec_2} we study in detail this scheme and provide its intersection numbers. Also a quotient scheme is considered, which turns out to arise from the point-graph of the generalized quadrangle. In Section \ref{sec_3}  we prove that any scheme having the same parameters as  the scheme based on the set of flags of a generalized quadrangle is necessarily such a scheme. In  Section 4 we find all possible fusions of this scheme and give a full description for those of class 2 and 3. Finally in  Section 5 we give a full description and a characterization of a 4-class symmetric fusion.
	
	It is worth pointing out that  the adjacency algebra associated with  the association scheme on the set of flags of a generalized polygon was considered in  \cite{hig}. In particular, Higman  provides an alternative proof of the Feit-Higman Theorem \cite{fh} by finding irreducible representations of this algebra.
	
	\section{The scheme on flags of a $\mathrm{GQ}(s,t)$}\label{sec_2}
	Let $\cS=(\cP, \cL,\mathrm I)$  be a GQ of order $(s, t)$.  A {\em flag} of $\cS$  is any pair $(p,L)\in\cP\times \cL$ with $p\mI L$. We denote the set of all flags of $\cS$ by $\Omega$. On the set $\Omega$ we consider the following relations, which, together with the diagonal relation $R_0$, partition the set $\Omega^2$:
	\begin{itemize}
		\item[$R_1$:] $((p,L),(q,M))\in R_1$ if and only if $p=q$ and $L\neq M$.
		\item[$R_2$:] $((p,L),(q,M))\in R_2$ if and only if $L=M$ and $p\neq q$. 
		\item[$R_3$:] $((p,L),(q,M))\in R_3$ if and only if $p$ and $q$ are collinear on $L$. 
		\item[$R_4$:] $((p,L),(q,M))\in R_4$ if and only if $p$ and $q$ are collinear on $M$. 
		\item[$R_5$:] $((p,L),(q,M))\in R_5$ if and only if $p$ and $q$ are collinear but neither on $L$ nor on $M$.
		\item[$R_6$:] $((p,L),(q,M))\in R_6$ if and only if $L\cap M=\{r\},\text{\ with\ }r\neq p,q$. 
		\item[$R_7$:] $((p,L),(q,M))\in R_7$ if and only if $L\cap M=\emptyset$ \text{and} $p$ \text{and} $q$ \text{are not collinear}.
	\end{itemize}
	
	We are going to prove that $\cX=(\Omega,\cR)$, where $\cR=\{R_0,R_1,\ldots,R_7\}$, is an imprimitive noncommutative   association scheme.   Note that $R_{3^*}=R_4$, so the scheme is not symmetric.

	\begin{remark}\label{rem_1}
		{\em 
			Let $\Omega^D$ be the set of all flags in the dual quadrangle $\cS^D$. By applying the point-line duality for GQ to relations $R_1,\ldots,R_7$ on $\Omega$, we get relations $R_{1^D},\ldots,R_{7^D}$ on $\Omega^D$, and  if $\cX$ is an association scheme so is  $\cX^D=(\Omega^D,\{R_{i^D}\}_{i=0}^7)$.
			Let $\Delta$ denote the map that associates   $(p,L)\in \Omega$ with $(L,p) \in \Omega^D$, and set $(R_i)^\Delta=\{(L,p):(p,L)\in R_i\}$. Then,
			\begin{equation}\label{eq_1}
				\begin{array}{ccc}
					(R_1)^\Delta=R_{2^D}, & (R_2)^\Delta=R_{1^D}, & (R_3)^\Delta=R_{4^D},\\ 
					(R_4)^\Delta=R_{3^D}, &  (R_5)^\Delta=R_{6^D}, & (R_6)^\Delta=R_{5^D},\\
					&	(R_7)^\Delta=R_{7^D}. &
				\end{array}
			\end{equation}
		}
	\end{remark}
	
	We now show that all of the intersection numbers $p_{ij}^k$ are well defined. 
	
	\begin{lemma}\label{lemma_1}
		The valencies $\eta_k=p_{kk^*}^0$, with $k\in I^*$, are as follows:
		$\eta_1=t$, $\eta_2=s$, $\eta_3=\eta_4=st$, $\eta_5=st^2$, $\eta_6=s^2t$, $\eta_7=(st)^2$.
	\end{lemma}
	
	\begin{proof}
		We calculate $\eta_1, \eta_2, \eta_3, \eta_4, \eta_5,\eta_6$ directly, obtaining $\eta_7$ by subtraction, as $|\Omega|=\sum\limits_{i=0}^{7}{\eta_i}$. Since $\cS$ has order $(s,t)$, then $\eta_1=t$ and  $\eta_2=s$. Let $(p,L)$ be any element in $\Omega$. Since there are $t$ lines distinct from $L$ which are incident with $p$, and each of such lines is incident with $s$ points different from $p$, we get $\eta_3=st$. Similarly, we get $\eta_4=st$. To  compute $\eta_5$, for any line $N$ incident with $p$ and different from $L$ consider a point $q$ incident with $N$. Then every pair $(q,M) \in \Omega$, with $M\neq N$, is $5$-related to $(p,L)$. From axioms (GQ1) and (GQ2), we get $\eta_5=st^2$. The value of $\eta_6$ is obtained by applying the point-line duality to $\cS$.
	\end{proof}

For the intersection  numbers $p_{ij}^{k}$ the following formulas are known  \cite[pp. 21--23]{cp}:

\begin{equation}\label{eq_3}
	\eta_{k}p_{ij}^k=\eta_{i}p_{jk^*}^{i^*}=\eta_{j}p_{k^*i}^{j^*},
\end{equation}

\begin{equation}\label{eq_2}
	p_{ij}^{k}=p_{j^*i^*}^{k^*}.
\end{equation}

By taking into account this equations and the Remark \ref{rem_1}, on the set $T$ of all triplets $(k\ i \ j)$, $1\leq k,i,j \leq 7$, we may define the maps:

\[
I:(k \ i \ j ) \mapsto (i^{*} \ j \ k^{*} )
\]
\[
S:(k \ i \ j )\mapsto (k^{*} \ j^{*} \ i^{*} )
\]
\[
D:(k \ i \ j ) \mapsto (k^{\delta} \ i^{\delta}\ j^{\delta})
\]
where 
\[
\begin{array}{rcll}
* & :&  i \mapsto i  &\ \ \ \mathrm{for\ } i \neq3, 4\\
& & 3 \leftrightarrow 4\\ [.2in]
	\delta& :&  i \leftrightarrow i+1  &\ \ \ \mathrm{for\ } i =1, 3,5\\
	& &7 \mapsto 7
\end{array}.
\]

Under composition, $I$, $S$, $D$ generate a permutation group of order 12 acting on  $T$:

\begin{center}
	$G=\<I,S,D\>=\{id,\ I, \ I^{2} , \ S,\ IS, \ I^2S,\  D, \ ID, \ I^2D,\ SD,\ ISD, I^2SD\}$.
\end{center}
By making the group $G$ act on the triplet $(k\ i\ j)$, $1\leq k ,i, j \leq 7$, and using Eqs. (\ref{eq_3})--(\ref{eq_2}) and Remark \ref{rem_1},  in Table \ref{tab_3} we report the action of each element $g \in G$ on $(k\ i \ j)$ and the intersection number corresponding to the triplet $(k\ i \ j)^g$; to make the notation simpler, $i^{\delta*}$ means $(i^{\delta})^*$, and $p^{k^\delta}_{i^{\delta}j^{\delta}}=f(t,s)$  if $p^k_{ij}=f(s,t)$:
	
{\renewcommand\arraystretch{1.7}
\begin{longtable}{|c|c|c|c|c|c|c|c|c|c|c|c|}
	\hline
	\endfirsthead
	\multicolumn{3}{r}{\textit{$($Continue to next page$)$}}
	\endfoot\multicolumn{12}{l}{\textit{$($Continued from previous page$)$}}
	\endhead
	\hline
	\multicolumn{3}{c}{}\\
	\caption{\small{The action of the elements of $G$ on the triplet $(k \ i \ j)$}}\endlastfoot
	\hline\label{tab_3}
	$I$ & $(i^{*} \ j \ k^{*} )$& $\frac{\eta_k}{\eta_i}p^k_{ij}$\\
	\hline
	$I^{2}$&$ (j^{*} \  k^{*} \ i )$&$\frac{\eta_k}{\eta_j}p^k_{ij}$\\
	\hline
	$S $ &$(k^{*} \ j^{*} \ i^{*} )$& $p^k_{ij}$\\
	\hline
	$IS $ & $ (j \ i^{*} \ k )$ & $\frac{\eta_k}{\eta_j}p^k_{ij}$\\
	\hline
	$I^{2}S$ & $(i\ k \ j^{*})$ & $\frac{\eta_k}{\eta_i}p^k_{ij}$\\
	\hline
	$D$ &$(k^{\delta} \ i^{\delta} \ j^{\delta})$ & $f(t,s)$\\
	\hline
	$ID$ &  $(i^{\delta*} \ j^{\delta}\ k^{\delta *} )$ & $\frac{\eta_{k^{\delta}}}{\eta_{i^{\delta}}}f(t,s)$ \\ 
	\hline
	$I^{2}D$ & $ (j^{\delta  *}\ k^{\delta  *}\  i^{\delta} )$ & $\frac{\eta_{k^{\delta}}}{\eta_{j^{\delta}}}f(t,s)$ \\
	\hline
	$SD$ & $(k^{\delta  *} \ j^{\delta  *} \ i^{\delta  *})$ & $f(t,s)$\\
	\hline
	$ISD$ & $(j^{\delta}\ i^{\delta  *}\ k^{\delta} )$& $\frac{\eta_{k^{\delta}}}{\eta_{j^{\delta}}}f(t,s)$ \\
	\hline
	$I^2SD$ & $(i^{\delta}\ k^{\delta} \ j^{\delta  *} )$& $\frac{\eta_{k^{\delta}}}{\eta_{i^{\delta}}}f(t,s)$ \\
	\hline
\end{longtable}
}

Therefore, it sufficies to compute 44 of the $7^3=343$  intersection numbers $p^k_{ij}$ of $\cX$, with $k,i,j\neq0$. These 44 intersection numbers, together with their orbit under $G$ are reported in Table \ref{tab_2}, where the $((k\ i \ j),g)$-entry, $g \in G$, is the triplet $(k\ i \ j)^g$; the  cell is left empty if the corresponding triplet has been previously found.

{\renewcommand\arraystretch{1.5}
	\scriptsize
	\begin{longtable}{|c|c|c|c|c|c|c|c|c|c|c|c|}
		\hline
		\endfirsthead
		\multicolumn{12}{r}{\textit{$($Continue to next page$)$}}
		\endfoot\multicolumn{12}{l}{\textit{$($Continued from previous page$)$}}
		\endhead
		\hline
		\multicolumn{12}{c}{}\\
		\caption{The orbits under the action of $G$ on the triplets $(k\ i\ j)$}\endlastfoot
		\hline\label{tab_2}
		$(k \ i \ j) $   & $I$       & $I^2$ & $S$       & $IS$      & $I^2S$ & $D$       & $ID$      & $I^2D$ & $SD$      & $ISD$     & $I^2SD$ \\
		\hline
		(1 1 1) &         &                      &         &         &                       & (2 2 2) &         &                       &         &         &                        \\ 
		\hline
		(1 1 2) & (1 2 1) & (2 1 1)              &         &         &                       & (2 2 1) & (2 1 2) & (1 2 2)               &         &         &                        \\ 
		\hline
		(1 1 3) & (1 3 1) & (4 1 1)              & (1 4 1) & (3 1 1) & (1 1 4)               & (2 2 4) & (2 4 2) & (3 2 2)               & (2 3 2) & (4 2 2) & (2 2 3)                \\ 
		\hline
		(1 1 5) & (1 5 1) & (5 1 1)              &         &         &                       & (2 2 6) & (2 6 2) & (6 2 2)               &         &         &                        \\ 
		\hline
		(1 1 6) & (1 6 1) & (6 1 1)              &         &         &                       & (2 2 5) & (2 5 2) & (5 2 2)               &         &         &                        \\ 
		\hline
		(1 1 7) & (1 7 1) & (7 1 1)              &         &         &                       & (2 2 7) & (2 7 2) & (7 2 2)               &         &         &                        \\ 
		\hline
		(1 2 3) & (2 3 1) & (4 1 2)              & (1 4 2) & (3 2 1) & (2 1 4)               &         &         &                       &         &         &                        \\ 
		\hline
		(1 2 4) & (2 4 1) & (3 1 2)              & (1 3 2) & (4 2 1) & (2 1 3)               &         &         &                       &         &         &                        \\ 
		\hline
		(1 2 5) & (2 5 1) & (5 1 2)              & (1 5 2) & (5 2 1) & (2 1 5)               & (2 1 6) & (1 6 2) & (6 2 1)               & (2 6 1) & (6 1 2) & (1 2 6)                \\ 
		\hline
		(1 2 7) & (2 7 1) & (7 1 2)              & (1 7 2) & (7 2 1) & (2 1 7)               &         &         &                       &         &         &                        \\ 
		\hline
		(1 3 3) & (4 3 1) & (4 1 3)              & (1 4 4) & (3 4 1) & (3 1 4)               & (2 4 4) & (3 4 2) & (3 2 4)               & (2 3 3) & (4 3 2) & (4 2 3)                \\
		\hline
		(1 3 4) & (4 4 1) & (3 1 3)              &         &         &                       & (2 4 3) & (3 3 2) & (4 2 4)               &         &         &                        \\ 
		\hline
		(1 3 5) & (4 5 1) & (5 1 3)              & (1 5 4) & (5 4 1) & (3 1 5)               & (2 4 6) & (3 6 2) & (6 2 4)               & (2 6 3) & (6 3 2) & (4 2 6)                \\ 
		\hline
		(1 3 6) & (4 6 1) & (6 1 3)              & (1 6 4) & (6 4 1) & (3 1 6)               & (2 4 5) & (3 5 2) & (5 2 4)               & (2 5 3) & (5 3 2) & (4 2 5)                \\ 
		\hline
		$(k \ i \ j) $   & $I$       & $I^2$ & $S$       & $IS$      & $I^2S$ & $D$       & $ID$      & $I^2D$ & $SD$      & $ISD$     & $I^2SD$ \\
		\hline
		(1 3 7) & (4 7 1) & (7 1 3)              & (1 7 4) & (7 4 1) & (3 1 7)               & (2 4 7) & (3 7 2) & (7 2 4)               & (2 7 3) & (7 3 2) & (4 2 7)                \\ 
		\hline
		(1 4 3) & (3 3 1) & (4 1 4)              &         &         &                       & (2 3 4) & (4 4 2) & (3 2 3)               &         &         &                        \\ 
		\hline
		(1 4 5) & (3 5 1) & (5 1 4)              & (1 5 3) & (5 3 1) & (4 1 5)               & (2 3 6) & (4 6 2) & (6 2 3)               & (2 6 4) & (6 4 2) & (3 2 6)                \\ 
		\hline
		(1 4 6) & (3 6 1) & (6 1 4)              & (1 6 3) & (6 3 1) & (4 1 6)               & (2 3 5) & (4 5 2) & (5 2 3)               & (2 5 4) & (5 4 2) & (3 2 5)                \\ 
		\hline
		(1 4 7) & (3 7 1) & (7 1 4)              & (1 7 3) & (7 3 1) & (4 1 7)               & (2 3 7) & (4 7 2) & (7 2 3)               & (2 7 4) & (7 4 2) & (3 2 7)                \\ 
		\hline
		(1 5 5) & (5 5 1) & (5 1 5)              &         &         &                       & (2 6 6) & (6 6 2) & (6 2 6)               &         &         &                        \\ 
		\hline
		(1 5 6) & (5 6 1) & (6 1 5)              & (1 6 5) & (6 5 1) & (5 1 6)               & (2 6 5) & (6 5 2) & (5 2 6)               & (2 5 6) & (5 6 2) & (6 2 5)                \\
		\hline
		(1 5 7) & (5 7 1) & (7 1 5)              & (1 7 5) & (7 5 1) & (5 1 7)               & (2 6 7) & (6 7 2) & (7 2 6)               & (2 7 6) & (7 6 2) & (6 2 7)                \\ 
		\hline
		(1 6 6) & (6 6 1) & (6 1 6)              &         &         &                       & (2 5 5) & (5 5 2) & (5 2 5)               &         &         &                        \\ 
		\hline
		(1 6 7) & (6 7 1) & (7 1 6)              & (1 7 6) & (7 6 1) & (6 1 7)               & (2 5 7) & (5 7 2) & (7 2 5)               & (2 7 5) & (7 5 2) & (5 2 7)                \\ 
		\hline
		(1 7 7) & (7 7 1) & (7 1 7)              &         &         &                       & (2 7 7) & (7 7 2) & (7 2 7)               &         &         &                        \\ 
		\hline
		(3 3 3) & (4 3 4) & (4 4 3)              &         &         &                       & (4 4 4) & (3 4 3) & (3 3 4)               &         &         &                        \\ 
		\hline
		(3 3 5) & (4 5 4) & (5 4 3)              &         &         &                       & (4 4 6) & (3 6 3) & (6 3 4)               &         &         &                        \\ 
		\hline
		(3 3 6) & (4 6 4) & (6 4 3)              &         &         &                       & (4 4 5) & (3 5 3) & (5 3 4)               &         &         &                        \\ 
		\hline
		(3 3 7) & (4 7 4) & (7 4 3)              &         &         &                       & (4 4 7) & (3 7 3) & (7 3 4)               &         &         &                        \\ 
		\hline
		(3 4 4) &         &                      & (4 3 3) &         &                       &         &         &                       &         &         &                        \\ 
		\hline
		(3 4 5) & (3 5 4) & (5 4 4)              & (4 5 3) & (5 3 3) & (4 3 5)               & (4 3 6) & (4 6 3) & (6 3 3)               & (3 6 4) & (6 4 4) & (3 4 6)                \\ 
		\hline
		(3 4 7) & (3 7 4) & (7 4 4)              & (4 7 3) & (7 3 3) & (4 3 7)               &         &         &                       &         &         &                        \\ 
		\hline
		(3 5 5) & (5 5 4) & (5 4 5)              & (4 5 5) & (5 5 3) & (5 3 5)               & (4 6 6) & (6 6 3) & (6 3 6)               & (3 6 6) & (6 6 4) & (6 4 6)                \\ 
		\hline
		(3 5 6) & (5 6 4) & (6 4 5)              & (4 6 5) & (6 5 3) & (5 3 6)               &         &         &                       &         &         &                        \\ 
		\hline
		(3 5 7) & (5 7 4) & (7 4 5)              & (4 7 5) & (7 5 3) & (5 3 7)               & (4 6 7) & (6 7 3) & (7 3 6)               & (3 7 6) & (7 6 4) & (6 4 7)                \\ 
		\hline
		(3 6 5) & (6 5 4) & (5 4 6)              & (4 5 6) & (5 6 3) & (6 3 5)               &         &         &                       &         &         &                        \\ 
		\hline
		(3 6 7) & (6 7 4) & (7 4 6)              & (4 7 6) & (7 6 3) & (6 3 7)               & (4 5 7) & (5 7 3) & (7 3 5)               & (3 7 5) & (7 5 4) & (5 4 7)                \\ 
		\hline
		(3 7 7) & (7 7 4) & (7 4 7)              & (4 7 7) & (7 7 3) & (7 3 7)               &         &         &                       &         &         &                        \\ 
		\hline
		(5 5 5) &         &                      & (6 6 6) &         &                       &         &         &                       &         &         &                        \\ 
		\hline
		(5 5 6) & (5 6 5) & (6 5 5)              &         &         &                       & (6 6 5) & (6 5 6) & (5 6 6)               &         &         &                        \\ 
		\hline
		(5 5 7) & (5 7 5) & (7 5 5)              &         &         &                       & (6 6 7) & (6 7 6) & (7 6 6)               &         &         &                        \\ 
		\hline
		(5 6 7) & (6 7 5) & (7 5 6)              & (5 7 6) & (7 6 5) & (6 5 7)               &         &         &                       &         &         &                        \\ 
		\hline
		(5 7 7) & (7 7 5) & (7 5 7)              &         &         &                       & (6 7 7) & (7 7 6) & (7 6 7)               &         &         &                        \\ 
		\hline
		(7 7 7) &         &                      &         &         &                       &         &         &                       &         &         &                        \\ 
		\hline
	\end{longtable}
}
Finally, by using the relation: 
\begin{equation}\label{eq_4}
	\sum\limits_{j=0}^7 p_{ij}^{k}=\eta_{i},
\end{equation}
we only need 28 intersection numbers; these are $p^1_{11}$, $p^1_{12}$, $p^1_{13}$, $p^1_{15}$, $p^1_{16}$, $p^1_{23}$, $p^1_{24}$, $p^1_{25}$,  $p^1_{33}$, $p^1_{34}$, $p^1_{35}$, $p^1_{36}$, $p^1_{43}$, $p^1_{45}$, $p^1_{46}$, $p^1_{55}$, $p^1_{56}$, $p^1_{66}$, $p^3_{33}$, $p^3_{35}$, $p^3_{36}$, $p^3_{44}$, $p^3_{45}$, $p^3_{55}$, $p^3_{56}$, $p^3_{65}$, $p^5_{55}$ and $p^5_{56}$ whose values are given in the following result.

\begin{proposition}\label{PropIn}
	The previous $28$  intersection numbers are all zeros except for $p^1_{11}=t-1$, $p^1_{23}=s$, $p^1_{35}=st$, $p^1_{43}=s(t-1)$, $p^1_{55}=st(t-1)$, $p^3_{35}=t(s-1)$, $p^3_{56}=st$, $p^5_{55}=t(s-1)$ and $p^5_{56}=s(t-1)$
\end{proposition}

\begin{proof}
	For any pair $((p,L),(p,M)) \in R_1$, we count the number of pairs $(z,N) \in \Omega$ such that $((p,L),(z,N)) \in R_{i}$ and $((z,N),(q,M))$ $ \in R_{j}$.\\
	Assume $i=j=1$. Then $z=p $ and $L \neq N  \neq  M$. From the axiom (GQ1) we get $p_{11}^1=t-1$. It is easy to see that $p^1_{1j}=0$ for $j\neq 1$.\\
	Assume $i=2$ and $j=3$. Then $z\neq p$, $L=N \neq M$ and $z$ and $p$ are collinear on $N$. From the axiom (GQ2) we get $p^1_{23}=s$. It is easy to see that $p^1_{24}=p^1_{25}=0$.\\
	Assume $i=3$ and $j=5$. Since $N\cap M = \emptyset$ and $p$ and $z$ are  collinear on $L$, from the axioms of a GQ, we get $p_{35}^1=st$. We also get $p^1_{33}=p^1_{34}=p^1_{36}=0$.\\	
	Assume $i=4$ and $j=3$. This implies that $z\neq p$, $N\cap L=\{p\}$ and $N \neq M$ (otherwise $((z,N),(p,M)) \in R_{2}$). From the axiom (GQ1) and (GQ2), we get $p^1_{43}=s(t-1)$. We also get $p^1_{45}=p^1_{46}=0$.\\
	Assume $i=j=5$. This implies that $p$ and $z$ are two distinct points collinear with a line, say $M'$, different from $L$, $N$ and also with $M$ (otherwise $((z,N),(p,M)) \in R_{4}$). From the axiom (GQ1) and (GQ2), we get $p^1_{55}=st(t-1)$. We also get $p^1_{56}=0$.\\
	By axiom (GQ3), $p^1_{66}=0.$\\
	
	For any pair $((p,L),(q,M)) \in R_3$, we count the number of pairs $(z,N) \in \Omega$ such that $((p,L),(z,N)) \in R_{i}$ and $((z,N),(q,M))$ $ \in R_{j}$.\\
	Assume $i=3$ and $j=5$. Then $p$, $z$ and $q$ are collinear on $L$, $L\cap N=\{z\}$, $L\cap M=\{q\}$ and $z\neq q$ (otherwise $((z,N),(q,M))$ $ \in R_{1}$). From the axiom (GQ1) and (GQ2), we get $p^3_{35}=t(s-1)$. We also get $p^3_{33}=p^3_{36}=0.$\\
	From the axiom (GQ3), $p^3_{44}=p^3_{45}=p^3_{55}=p^3_{65}=0.$\\
	Assume $i=5$ and $j=6$. Then $z$ and $p$ are collinear on $M'$ with $M' \neq  L,  M$, and $M \cap N=\{r\}$, with $r\neq z,q$. From the axioms of a generalized quadrangle, we get $p^3_{56}=st$.\\

	Finally, for any pair $((p,L),(q,M)) \in R_5$, with $p$ and $q$ collinear on a line $L'$ different from both $L$ and $M$, we count the number of pairs $(z,N) \in \Omega$ such that $((p,L),(z,N)) \in R_{5}$, with $p$ and $z$ are collinear on $M'$ different from $L$ and $N$, and $((z,N),(q,M))$ $ \in R_{j}$.\\
	Assume $j=5$. Then $L'=M'$. If $N=L'$, then  $((p,L),(z,L')) \in R_4$. So $N\neq L'$. From the axiom of a (GQ1), we get $p^5_{55} = t(s-1)$. \\
	Assume $j=6$. Then $L'\neq M'$. Fix any line $M'$ incident with $p$ and different from $L$ and $L'$. For each point $r $ incident with $M$, $r \neq q$, there is a unique flag $(z,N)$ such that $r\, {\rm I}\, N\, {\rm I}\, z\, {\rm I}\,  M'$. So, for the given line $M'$, there are $s$ flags $(z,N)$ $6$-related  to $(q,M)$; hence $p^5_{56}=s(t-1)$. 
\end{proof}

	\begin{corollary}\label{lemma22}
		The intersection numbers of the association scheme $\cX $  are collected in the following matrices $L_k$ whose $(i,j)-$entry is $p_{ij}^k$:
		\[
		L_1=\begin{pmatrix}
			0 & 	1 		 &	   0   &			0  		 & 	0		 &		 	0	 &0	&0	  \\
			1  &    t-1   &     0    &  0   & 0   & 0  &0 &0  \\
			0  &  0		&     0   &   s   & 0  &  0 & 0  & 0 \\
			0  &  0		&  	  0    &     0 &  0 &  st & 0  &0   \\
			0  &  0		&     s    &   s(t-1) &   0 &  0 & 0  & 0 \\
			0  &  0		&  	  0	   &   0  &   st & st(t-1)  & 0&  0   \\
			0  &  0		&     0	   &     0 &   0&  0 &0& s^{2}t  \\
			0  &  0		&     0	   &     0&  0 & 0 & s^{2}t   & s^{2}t(t-1)  \\
		\end{pmatrix}
		\]
		\medskip
		\[
		L_2=\begin{pmatrix}
			0	& 		0		 &		1	   &		 	0	  &		  0		 & 		 0		 &		 		 0		 &		0   \\
			0   &     0 &    0     &    0 & t   & 0   & 0  & 0   \\
			1  &  0 &  s-1 & 0  & 0 & 0 & 0  & 0 \\
			0 & t  &  0 & 0 &  t(s-1)  & 0  & 0  & 0  \\
			0  &  0 &  0 & 0 & 0  & 0  & st  & 0  \\
			0 &  0 &  0 & 0 &  0 &  0 &  0 & st^{2}  \\
			0 &  0 &  0 & st  &   0 &0 &  st(s-1)  &0  \\
			0& 0 &  0 & 0  & 0  & st^{2} & 0  & st^{2}  (s-1) \\
		\end{pmatrix}
		\]
		\medskip
		\[
		L_3=\begin{pmatrix}
			0 &	   0	 	&0	   		&1  		 &0  		     &0 			    &0		 		 		&0					 \\
			0 &    0    	& 0    		& 0	   		  & 0   		 &t 		       &0                       &0  	 				 \\
			0 & 	1 	&	0   	 & s-1   		  & 0      	&0		          &0            		    &0					 \\
			1 &	 t-1 	& 0 		 & 0   		 &  0  	   		 &t(s-1)  	   & 0		      & 0 					  \\
			0 &  	0		& 	0  		 &  0 	     &0         	 &0			       & 0					    &st				  	   	 \\
			0 & 	0  		&0   		& 0    	   	 &0     	  	 &0               & st  				    &st(t-1)   		  	 \\
			0 &		0   	&s 	 		 & s(t-1)    	 & 0	   &0          	     & 0  				      &st(s-1)				  \\
			0 & 	0  		&0 	 		& 0  		 & st  	 		 &st(t-1)       &  st(s-1) 			   & st(s-1)(t-1)		 \\
		\end{pmatrix}
	\]
	\medskip
		\[
		L_4=\begin{pmatrix}
			0	& 			0	 &		0	   &		 0		  &		1  		 & 		 	0	 &		 0		 		 &	0	  \\
			0 &    0   &    1   &0   &   t-1  & 0   & 0  &  0   \\
			0  & 0  &0  & 0& 0& 0    &   s&0 \\
			0 & 0  & 0  & 0  & 0  & 0  & 0  & st  \\
			1 &0   & s-1  &  0&0 &0   &  s(t-1)    & 0  \\
			0 & t  & 0  &0 & t(s-1)   & 0  & 0  &st(t-1)   \\
			0  & 0  & 0  & 0  & 0  & st  & 0  &st(s-1)  \\
			0 & 0  & 0  &  st  &0 & st(t-1) &  st(s-1) & st(s-1)(t-1)  \\
		\end{pmatrix}
		\]
		\medskip
		\[
				L_5=\left(\begin{array}{cccccccc}
			0 &	   0	 	&0	   	  &0  	    	 &0  		      &1 			     &0		 		 		&0					 \\
			0 &    0    	& 0    		& 1	   		  & 0   		  &t-1 		          &0                     &0  	 				 \\
			0 & 	0 	 	&	0   	 & 0   		  & 0        	  &0		        &0            		    &s					 \\
			
			0  &   	 0	    & 0 		 & 0   		   &  0  	   		 &0 		  	    & s		                &s(t-1)	 					  \\
			0&  	1		& 	0  		 &  s-1 	   &0               &0			        & 0					    &s(t-1)				  	   	 \\

			1 & 	t-1  	&0   		& 0    	   	 &0         	 &	 t(s-1)        &  s(t-1)		      &s(t-1)^2   		  	 \\
			0 &		0   	&0 	        & 0  	  	  & s	 		 &  s(t-1)          & s(s-1)	  			 &s(s-1)(t-1)				  \\
			0 & 	0  		&s 	 		& s(t-1)	& s(t-1)	 &s(t-1)^2     & s(s-1)(t-1)		 & s(s-1)(t^2-t+1)		 \\
		\end{array}\right)
		\]
		\medskip
		\[
		L_6=\left(\begin{array}{cccccccc}
			0 &	   0	 	&0	   	  &0  		 &0  		      &0 			     &1		 		 		&0					 \\
			0 &    0    	& 0    		& 0	   		  & 0   		&0 		         &0                       &t  	 				 \\
			0 & 	0 	 	&	0   	 & 0   		  & 1        	   &0		         &s-1            		    &0					 \\
			0  &   	 0	  & 1 		  & 0   		 &  t-1  	   		 &0 		  	   & 0		              & t(s-1)	 					  \\
			0&  	0		& 0  		 &  0 	     &0         	 &t			       & 0					    &t(s-1)				  	   	 \\
			0 & 	0  		&0   		& t    	   	 &0     	  	 &t(t-1)              &t(s-1)  		    &t(s-1)(t-1)   		  	 \\
			1 &		0   	&s-1 	 & 0    	  & 0	 		    &t(s-1)            &  s(t-1)			  &t(s-1)^{2}				  \\
			0 & 	t  		&0 	 		& t(s-1)  	& t(s-1)    &   t(s-1)(t-1)  & t(s-1)^{2}     & t(t-1)(s^2-s+1)		 \\
		\end{array}\right)
	\]
	\medskip
		\[
		{\scriptsize
			\hspace{-1.5cm}L_7=\left(\begin{array}{cccccccc}
				0 &	   0	 	&0	   	  &0  	    	 &0  		      &0 			     &0		 		 		&1					 \\
				0 &    0    	& 0    		& 0	   		  & 0   		  &0 		          &1                    &t-1  	 				 \\
				0 & 	0 	 	&	0   	 & 0   		  & 0        	  &1		         &0            		    &s-1					 \\
				
				0  &   	 0	    & 0 		 & 1   		 &  0  	   		 &t-1 		  	    &s-1 		             &(s-1)(t-1)			  \\
				0&  	0		& 	0  		 &  0 	      &1         	  &t-1			    &s-1					 &(s-1)(t-1)			\\
				
				0 & 	0	  	&1   		& t-1    	  &t-1         	 &	(t-1)^2   & (s-1)(t-1)		  &(s-1)(t^2-t+1)		 \\
				0 &		1   	&0 	        & s-1  	  	  & s-1	 		&  (s-1)(t-1)      & (s-1)^2	    &(s^2-s+1)(t-1)	    \\
				
				1 & 	t-1     &s-1 	 & (s-1)(t-1)	&(s-1)(t-1)	 & (s-1)(t^2-t+1) & (s^2-s+1)(t-1) &1 - s + s^2 - t - s^2 t + t^2 - s t^2 + s^2 t^2 	 \\
			\end{array}\right)
			}	
		\]
	\end{corollary}

\begin{proof}
	Each $p^k_{ij}$ is computed by using Proposition  \ref{PropIn} and Tables \ref{tab_3} and \ref{tab_2}.
\end{proof}

	\begin{theorem}\label{th_1}
		The pair $\cX=(\Omega,\{R_i\}_{i=0}^{7})$ is a noncommutative, imprimitive association scheme of order $(s+1)(t+1)(st+1)$ and  class $7$. The intersection numbers of this scheme are polynomials in $s$ and $t$. 
	\end{theorem}
	\begin{proof}
		From Lemma \ref{lemma22},  $\cX$ is an association scheme of order $(s+1)(t+1)(st+1)$ and  class 7. Since $p^1_{45}\neq p^1_{54}$, $\cX$ is not commutative.
		
		For every fixed $x \in \Omega$, the set $R_1(x)\cup\{x\}$ is the vertex set of a  $\{0,1\}$-clique of size $t+1$. Since $p^1_{11}=t-1$, such a clique is maximal, and the  basis graph $(\Omega,R_1)$ is the disjoint union of $(st+1)(s+1)$ maximal $\{0,1\}-$cliques. This implies  that $R_0\cup R_1$ is a non-trivial parabolic of $\cX$, whose equivalence classes are the maximal $\{0,1\}-$cliques. Since $p^2_{22}=s-1$, the same holds for $R_0\cup R_2$. Hence, $\cX$ is imprimitive.
	\end{proof}
	
	\begin{remark}\label{remark_2}
	{\em If $s=t=1$, then $\cX=(\Omega,\{R_i\}_{i=0}^{7})$ is a thin association scheme with
			\[
			\Omega=\{ (a,A), (b,A), (b,B), (c,B), (c,C), (d,C), (d,D), (a,D)  \}.
			\]
			Direct computation shows that, with respect to the complex product, $R_1$ and $R_3$ have order $2$ and $4$, respectively, and $R_1R_3R_1=R^{-1}_3=R_4$. This yields that $\cR=\{R_i\}_{i=0}^{7}$, endowed with the complex product, is isomorphic to the dihedral group $D_{8}$.
		}
	\end{remark}
	
	From now on,  ``$\{0, i\}$-clique'' will stand for  ``maximal $\{0, i\}$-clique''. In addition,  to make the notation lighter, we identify a clique $C$ with its vertices; for any clique $C$, we will write $x\in C$ to denote a vertex $x$ of $C$, if no confusion arises. 
	
	Since $e_1=R_0\cup R_1$ and $e_2=R_0\cup R_2$ are non-trivial parabolics of $\cX$,  it is possible to  construct the quotient  schemes on $\Omega/e_1$ and on $\Omega/e_2$ by following  \cite[Section 3.1.2]{cp}. It is evident that the elements of $\Omega/e_1$ are the $\{0, 1\}$-cliques and those of $\Omega/e_2$ are the $\{0, 2\}$-cliques.
	
	We call the elements of $\Omega/e_1$ \emph{point-cliques} and those of $\Omega/e_2$ \emph{line-cliques}.  We say that $C_1 \in \Omega/e_1$ and $C_2 \in \Omega/e_2$ are \emph{incident}, and we will write $C_1$ ${\rm I} $ $C_2$, if $C_1 \cap C_2 \neq \emptyset$. Note that, in this case, $|C_1 \cap C_2|=1 $, as $\cR$ is a partition of $\Omega \times \Omega$. In addition,  every $x\in \Omega$ is contained in a unique point-clique and a unique line-clique, which will be denoted by  $C_1(x)$ and  $C_2(x)$, respectively.

	\begin{lemma}\label{lemma11}
		Let $(x,y) \in R_2$. Then $(w,z) \in R_2\cup R_3\cup R_4\cup R_5$ for every $(w,z) \in C_1(x)\times C_1(y)$.
	\end{lemma}
	\begin{proof}
		Let $w \in C_1(x) \setminus \{x\}$ and $z \in C_1(y) \setminus \{y\}$. 
		Since $p^2_{1i}\neq 0$ only for $i=4$, then $(w,y) \in R_4$ and $(x,z) \in R_3$ (as $p^2_{i1}\neq 0$ only for $i=3$). Since $p^4_{i1}\neq0$ only for $i=5$, then $(w,z)\in R_5$  (and the same holds if we consider $p^3_{1i}$). 
	\end{proof}
	
	\begin{lemma}\label{lemma12}
		Let $(x,y) \in R_6$. Then $(w,z) \in R_6\cup R_7$ for every $(w,z) \in C_1(x)\times C_1(y).$
	\end{lemma}
	\begin{proof}
		Let $w \in C_1(x) \setminus \{x\}$ and $z \in C_1(y) \setminus \{y\}$.
		Since $p^6_{1i} (= p^6_{i1})\neq 0$ only for $i=7$, then $(w,y) \in R_7$ and $(x,z) \in R_7$. Since $p^7_{i1}(= p^7_{1i})\neq0$ for $i=6,7$, then $(w,z)\in R_6 \cup R_7$. 
	\end{proof}
	
	By the previous lemmas, we define the following nontrivial relations on $\Omega/e_1$:
	\begin{itemize}
		\item[$\bar R_1$:] $(C_1,C'_1)\in \bar R_1$ if and only if $C_1\times C'_1 \subseteq R_2\cup R_3\cup R_4\cup R_5  $.
		\item[$\bar R_2$:] $(C_1,C'_1)\in \bar R_2$ if and only if $C_1\times C'_1 \subseteq R_6\cup R_7$.
	\end{itemize}
	\begin{proposition}
		The  basis graph $(\Omega/e_1,\bar R_{1})$ of the quotient scheme $\cX_{\Omega/e_1}$ is the point-graph of the generalized quadrangle $\cS$.
	\end{proposition}
	\begin{proof}
		Let $(C_1,C'_1)\in\bar R_1$, and  $(x,y)\in (C_1\times C'_1)\cap R_2$. This implies  
		\[
		C_2(x)=\{z\in\Omega:(x,z)\in R_2\}=C_2(y).
		\]
		Therefore, $C_1\,{\rm I'}\,  C_2(x)\,{\rm I'}\, C'_1$, i.e. $C_1$ and $C'_1$ are two collinear point-cliques.
		
		On the other hand, it is easily seen that if $C_1$ and $C'_1$ are two point-cliques that are both incident with  a line-clique $C_2$, then $(C_1,C'_1)\in \bar R_1$.  Therefore $(\Omega/e_1,\bar R_{1})$  is the point-graph of the generalized quadrangle $\cS$.
	\end{proof}
	
	By applying very similar arguments, it is possible to prove that the quotient scheme $\cX_{\Omega/e_2}$ is the point-graph of the generalized quadrangle $\cS^D$.

	
	\section{Reconstructing the generalized quadrangle from the scheme $\cX$}\label{sec_3}
	
	Let $\cX'=(\Omega',\{R'_i\}_{i=0}^{7})$ be an association scheme that is algebraically isomorphic to $\cX=(\Omega,\{R_i\}_{i=0}^{7})$ via the isomorphism $\phi$. To make notation simpler, we set $R'_{i}=\phi(R_i)$, for $i=0,\ldots,7$. By \cite[Prop. 2.3.25]{cp} $e'_{1}=R'_{0}\cup R'_{1}$ and $e'_{2}=R'_{0}\cup R'_2$ are parabolics of $\cX'$.  
	
	Our aim is to reconstruct the generalized quadrangle with parameters $(s,t)$ from $\cX'$.  Set $\cP'=\Omega'/e'_1 $ and $\cL'=\Omega'/e'_2 $.  We say that $C_1 \in \cP'$ and $C_2 \in \cL'$ are \emph{incident}, and we will write $C_1$ ${\rm I'} $ $C_2$, if $C_1 \cap C_2 \neq \emptyset$, that is $|C_1 \cap C_2|=1 $.
	
	We are going to show that  $\cS'=(\cP',\cL',  {\rm I'})$ is a generalized quadrangle of order $(s,t)$. 
	
	Since every point-clique has $t+1$ vertices, each of which is on a unique line-clique, it follows that every point-clique is incident with $t+1$ line-cliques. Similarly, we find that every line-clique is incident with $s+1$ point-cliques. So, from the maximality of the cliques, axioms (GQ1) and (GQ2) are satisfied. 
	
	\begin{lemma}\label{lemma9}
		Let $C_1 \in \cP'$  and $C_2 \in \cL'$ be incident, with common vertex $z$. Then $(x,y) \in R'_4$ for all $x \in C_1\setminus \{z\}$ and $y \in C_2\setminus \{z\}$. Conversely,  if $(x,y) \in R'_4$ then $C_1(x)$ and $ C_2(y)$ are incident.
		
	\end{lemma}
	\begin{proof}
		Since $p^2_{1i}\neq0 $ only for $i=4$ the first part of the statement follows. Conversely, let  $(x,y) \in R_4^{\prime}$. Then $p^4_{i2}\neq 0$ only for $i=1,4.$ In particular, $p^4_{12}=1$ implies that  $| C_1(x) \cap C_2(y)|=1$, i.e. $C_1(x)$ ${\rm I'}$ $C_2(y)$.
	\end{proof}
	
	\begin{lemma}\label{lemma10}
		Let $C_1 \in \cP'$ and $C_2 \in \cL'$  be not incident. Then, there exists at most one pair $(x,y) \in C_1\times C_2$ such that $(x,y)\in R_3^{\prime}$.
	\end{lemma}
	\begin{proof}
		Assume that there are two distinct pairs  $(x,y), (x',y') \in C_1\times C_2$ both in $R_3^{\prime}$. If $x\neq x'$ and $y=y'$ we should have $p^1_{34}\neq 0$; a contradiction. Similarly if $y \neq y'$ and $x=x'$ we get a contradiction since $p^2_{43}= 0$. Let $x\neq x'$ and $y \neq y'$. Since $p^1_{3j} \neq 0$ for $j=5$ then $(y,x') \in R_5^{\prime}$, giving $p^2_{53}\neq 0$; a contradiction. 
	\end{proof}
	
	\begin{theorem}
		Let $\cX'=(\Omega',\{R'_i\}_{i=0}^{7})$ be an association scheme that is algebraically isomorphic to $\cX=(\Omega,\{R_i\}_{i=0}^{7})$. Then $\cX'$ is the association scheme constructed on the flags of a generalized quadrangle.
	\end{theorem}
	\begin{proof}
		By keeping in mind the above notation, let $R'_{i}=\phi(R_i)$, for $i=0,\ldots,7$.\\
		We remark that if $(x,y) \in R_3^{\prime}$ then 
		\[
		C_1(x) \, {\rm I'} \, 	C_2(x) \, {\rm I'} \, 	C_1(y) \, {\rm I'} \, 	C_2(y),
		\]
		by Lemmas \ref{lemma9} and \ref{lemma10}. Fix $C_1 \in \cP'$ and $x \in C_1$, so that $C_1(x)=C_1$. By Lemma \ref{lemma9}, for every $y \in R_3^{\prime}(x)$ the clique $C_1$ and $C_2(y) $ are not incident. By Lemma \ref{lemma10}, the set $\cL_3^{\prime}(x)= \{C_2(y) : y \in R_3^{\prime}(x)\}$ consists of $\eta_3=st$ pairwise distinct line-cliques, and each of them is not incident with $C_1$ by Lemma \ref{lemma9}. Again by Lemma \ref{lemma10}, $\cL_3^{\prime}(x) \cap \cL_3^{\prime}(x')=\emptyset $, for $x,x' \in C_1$ with $x\neq x'$. It follows that there are $(t+1)st$ line-cliques in $\cX'$ which are not incident with $C_1$. 
		Set 
		\[
		\cL^{\prime}(C_1)=\bigcup_{x \in C_1} \{C_2(y) : y \in R_3^{\prime}(x)\} \cup \{C_2(x) : x \in C_1\}.
		\]
		
		Since $|\cL^{\prime}(C_1)|=st(t+1)+(t+1)=(st+1)(t+1)=|\cL^{\prime}|$, we get $\cL^{\prime}(C_1)=\cL^{\prime}$. This implies that for every clique $C_2$ not incident with $C_1$ there exist unique line-clique $D_2$ and point-clique $D_1$ such that 
		\[
		C_1 \, {\rm I'} \, 	D_2 \, {\rm I'} \, 	D_1 \, {\rm I'} \, 	C_2.
		\]
		From the arbitrariness of the choice of $C_1 \in \cP'$, the axiom (GQ3) holds in $\cS'=(\cP',\cL', {\rm I'})$.
	\end{proof}
	
	\section{The fusions of the scheme $\cX$}
	
	{
		Let $\cX=(X, \{R_i\}_{i=0}^d)$ be an association scheme, and  $\{\Lambda_0, \Lambda_1,\ldots, \Lambda_e\}$, $e \geq 2$ be a partition of $\{0,1,\ldots, d\}$ such that $\Lambda_0 = \{0\}$. Set $S_l=\bigcup_{l'\in\Lambda_l}{R_{l'}}$, $l=0,\ldots,e$. If $\cY=(X,\{S_l\}_{l=0}^e)$ is an association scheme, then $\cY$ is called a {\em non-trivial fusion} of the association scheme $\cX$.
		The partition $\{\Lambda_1,\ldots, \Lambda_e\}$ of $\{1,\ldots, d\}$ gives rise to a fusion of $\cX$ if the following criterion is satisfied: for $l,l',l''\in \{1,\ldots,e\}$, the equations
		\begin{equation}\label{eq_6}
			\sum_{{\substack{i\in\Lambda_{l'} \\ j\in\Lambda_{l''}}}}{p^k_{ij}}=\sum_{{\substack{i\in\Lambda_{l'} \\ 
						j\in\Lambda_{l''}}}}{p^{k'}_{ij}}, 
		\end{equation}
		hold for all  $k,k'\in \Lambda_l$.
		
		Let $\cX$ be the association scheme on the flags of a GQ$(s,t)$.		
		For $k,i,j \in \{1,\ldots,7\}$, let $f^k_{ij}(x,y)$ be the polynomial such that $p^k_{ij}=f^k_{ij}(s,t)$. Let $\{\Lambda_1,\ldots, \Lambda_e\}$ be a non-trivial fusion of $\cX$. Then Eqs. (\ref{eq_6}) can be written as \begin{equation}\label{eq_7}
			\sum_{{\substack{i\in\Lambda_{l'} \\ j\in\Lambda_{l''}}}}{f^k_{ij}(s,t)}=\sum_{{\substack{i\in\Lambda_{l'} \\ 
						j\in\Lambda_{l''}}}}{f^{k'}_{ij}(s,t)}.
		\end{equation}
		Under the point-line duality described in Remark \ref{rem_1}, $\{\Lambda^D_1,\ldots, \Lambda^D_e\}$, where $\Lambda^D_l=\{i^D : i \in \Lambda_l\}$, is a partition of $\{1,\ldots,7\}$. Recall that $\{1^D,\ldots,7^D\}$ are the relations of the association scheme on the  flags of the dual quadrangle $\cS^D$, with  $p^{k^D}_{i^Dj^D}=f^k_{ij}(t,s)$. Since Eqs. (\ref{eq_7}) still hold if we interchange $s$ with $t$, we get that  $\{\Lambda^D_1,\ldots, \Lambda^D_e\}$ gives a fusion of $\cX$ viewed as the scheme constructing on $\cS^D$.
		
		Let $\{\Lambda_1,\ldots, \Lambda_e\}$, $e \geq 2$, be a partition of $\{1,\ldots, 7\}$. Since $R_3$ and $R_4$ are the only non-symmetric basis relations on $\Omega$, it easy to see that either $\{3,4\}\subseteq \Lambda_i$ for some $i=1,\ldots,e$, or the singletons $\{3\}$ and $\{4\}$ are elements of the partition. 
		Taking into account this remark and equations  (\ref{eq_6}), we use the computer algebra system Mathematica \cite{math} to find all the partitions $\{\Lambda_1,\ldots, \Lambda_e\}$ of $\{1,\ldots,7\}$ such that the corresponding association scheme $\cY$ is a fusion of $\cX$. Since for $(s,t)=(1,1)$ the scheme is completely described in Remark \ref{remark_2}, we just consider fusions arising from GQ$(s,t)$ with $(s,t)\neq(1,1)$.
		
		All fusions of the scheme $\cX$, up to duality in the sense described above, are given in Table \ref{tab_1}. It turns out that there is no fusion such the corresponding partition contains the singletons $\{3\}$ and $\{4\}$.

		{
			\renewcommand\arraystretch{1.3}
			\begin{table}[!h]\hspace{-0.7cm}
				\scriptsize
				\begin{tabular}{|c|c|c|c|c|c}
					\hline
					Partition &\begin{tabular}[c]{@{}c@{}}  Feasible values\\ for s and t\end{tabular} &  Basis graph & Type & Reference \\
					\hline
					\{1,2,3,4,7\} \{5,6\}& $s=t=2$ & $(\Omega,\{1,2,3,4,7\})$  & $pg(4,6,3)$ &  \cite{bro} \\
					\hline		
					\{1,3,4,6\} \{2,5,7\}& $s \in \mathbb{N}, t=1$ &$(\Omega,\{2,5,7\})$ &$2K_{(s+1)^2} $   &\cite{vmb}\\
					\hline
					\{1,2,3,4,6,7\} \{5\}& $s\in \mathbb{N}, t=1$& $(\Omega,\{5\})$ & $2(s+1)K_{s+1}$ & \cite{vmb}  \\
					\hline
					\{1\} \{2,3,4,5,6,7\} & $s,t \in \mathbb{N}$& $(\Omega,\{1\})$  & $(st+1)(s+1)K_{t+1}$ & \cite{vmb} \\
					\hline
					\{1,2,7\} \{3,4\} \{5,6\} & $s=t=2$  & $(\Omega,\{5,6\})$ & Gerwitz$_2(x)$ & \cite[p. 93]{vd}\\
					\hline
					\{1\} \{2,5,6\} \{3,4,7\} & $s= 3, t=1$& $(\Omega,\{2,5,6\})$ &Halved 6-cube graph  &  \cite[p. 92]{vd}	\\				
					\hline
					\{1,6\} \{2,5,7\} \{3,4\}  & $s=3,  t=1$ &$(\Omega,\{3,4\})$  & Folded 6-cube  & \cite[p. 92]{vd},\cite{distreg}\\
					\hline
					\{1\} \{2,5,7\} \{3,4,6\} & $s \in \mathbb{N},\ t=1$& $(\Omega,\{3,4,6\})$ & $R(2,(s+1)^2)$  &\cite[p. 88]{vd}\\
					\hline
					\{1,3,4,6\} \{2,5\} \{7\}& $s\in \mathbb{N}, t=1$ & $(\Omega,\{2,5\})$ &  $H(2,s+1)$ &\cite[p. 88]{vd},  \cite{vmb}\\
					\hline
					\{1,3,4,6\} \{2,7\} \{5\}& $s\in \mathbb{N}, t=1$ &$(\Omega,\{1,3,4,6\})$ &  SRG$(2(s+1),s+1,0,s+1) \otimes J_{s+1}$   &\cite[p. 88]{vd}\\
					\hline
					\{1\} \{2,3,4,5\}  \{6,7\}& $s, t\in \mathbb{N}$& $(\Omega,\{2,3,4,5\})$  & SRG$((st+1)(s+1), s(t+1), s-1,t+1 ) \otimes J_{t+1}$ &\cite[p. 88]{vd}\\
					\hline
					\{1,2\} \{3,4\} \{5,6\} \{7\}& $s=t,  s \in \mathbb{N}$ & $(\Omega,\{1,2\})$ & \parbox[c][.3in][c]{2.2in}{Incidence graph  of the dual  of the double of $\cS$; see Section 5} & \\
					\hline
					\{1,3,4,6\} \{2\} \{5\} \{7\}&$s\in \mathbb{N}, t=1$ &  &  &\\
					\hline
					\{1\} \{2,5\} \{3,4\} \{6\} \{7\}& $s \in \mathbb{N},\ t=1$ & &  & \\
					\hline
				\end{tabular}	
				\caption{\small{The non-trivial fusions of the association scheme $\cX$}}\label{tab_1}
			\end{table}
		}
\begin{remark}
{\em In \cite{leo}, Leonard introduced systematic Gr\"obner basis methods to funding all fusions of the scheme  on flags of a generalized polygon. To compare Leonard's result on generalized quadrangles with ours, it is necessary to keep in mind that the flag adjacency matrices $A_1,A_2,A_5$ and $A_6$ in \cite{leo} are the matrices for our relations $R_2,R_1,R_6$ and $R_5$, respectively.}
\end{remark}

		\section{The fusion of $\cX$ from $\{\{1,2\},\{3,4\},\{5,6\},\{7\}\}$}
		
		In the following, for any given pair $(x,y)\in\Omega\times\Omega$, we set	$\cV^{(x,y)}_{ij}=\{z\in \Omega:(x,z)\in S_i,(z,y)\in S_j\}$. Obviously,  $|\cV^{(x,y)}_{ij}|=p^k_{ij}$ if $(x,y)\in S_k$.
		
		Let $ \tilde\cR=\{R_0$, $S_{1}=R_1\cup R_2$, $S_{2}=R_3\cup R_4$, $S_{3}=R_5\cup R_6$, $S_{4}=R_7\}$. In this section, we study in detail the fusion $\cY=(\Omega,\tilde\cR)$ and we will prove that this scheme is characterized by its parameters.

		\begin{theorem}\label{th_3}
			The pair $\cY=(\Omega,\tilde\cR)$ is a symmetric, primitive association scheme of order $(s+1)^2(s^2+1)$ and  class $4$. The valencies are 
			\[
			\eta_{1}=2s,\ \ \eta_{2}=2s^2,\ \ \eta_{3}=2s^3,\ \  \eta_{4}=s^4.
			\]
			The intersection numbers are collected in the following matrices whose $(i,j)-$entry is $p_{ij}^{k}$:
			\[
			M_{1}=\begin{pmatrix}
				0 & 	1 		 &	   0   &			0  		 & 	0		  \\
				1  &    s-1   &     s    &  0  & 0   \\
				0  &  s		&     s(s-1)   &   s^2   & 0  \\
				0  &  0		&  	  s^2   &     s^2(s-1)   &  s^3   \\
				0  &  0		&     0    &   s^3 &   s^3(s-1)  \\
			\end{pmatrix}
			\]
			\medskip
			\[
			M_{2}=\begin{pmatrix}
				0 &	   0	    		&1  		 &0  					    &0		 					 \\
				0 &    1    	& s-1   		& s	   		  & 0    				 \\
				1  &	 s-1 	& 0 				 &s(s-1)  	   & s^2		  \\
				0 &		s 	 	& s(s-1)    & s^2	      &2s^2(s-1)				  \\
				0 & 	0  		 & s^2  	 		 &2s^2(s-1)   	   & s^2(s-1)^2		 \\
			\end{pmatrix}
			\]	
			\medskip
			\[
			M_{3}=\left(\begin{array}{cccccccc}
				0 &	   0	 	&0	   	  &1 	    	 &0\\
				0 &    0    & 1	   	  &s-1 		          &s \\
				0 & 	1 	 	&	s-1   	 & s   		  & 2s(s-1)  \\
				1 & 	s-1  	&s   		& 4s(s-1)    	   	 &2s(s-1)^2   		  	 \\
				0 & 	s  		&2s(s-1)	&2s(s-1)^2     & s(s-1)(s^2-s+1)		 \\
			\end{array}\right)
			\]
			\medskip
			\[
			M_{4}=\left(\begin{array}{cccccccc}
				0 &	   0	 	&0	   	  &0  	    	 &1 \\
				0 &    0    	& 0    		& 2   		  & 2(s-1)\\
				0 & 	0 	 	&	2		         &4(s-1)        &2(s-1)^2					 \\
				0 & 2  		& 4(s-1)    	   &	4(s-1)^2   &2(s-1)(s^2-s+1)		 \\
				1 & 	2(s-1)     &2(s-1)^2 & 2(s-1)(s^2-s+1) & s^4-2s^3+2s^2-2s+1  \\
			\end{array}\right)
			\]
		\end{theorem}

		\begin{proof}
			A direct computation shows that the given partition and the corresponding intersection numbers of $\cX $ satisfy Eqs. (\ref{eq_6}). This implies that $\cY$ is a (nontrivial) fusion of $\cX$ and Eqs. (\ref{eq_6}) provides, at the same time,  the intersection numbers of it. 
			
			We now show that $\cY$ is a primitive scheme by checking that every basis graph $\Gamma_i=(\Omega,S_i)$, for $i=1,2,3,4$, is connected. By definition, $\Gamma_i$ is connected  if, for any two non-adjacent vertices $x$ and $y$, there exists an $S_i$-path  from  $x$ to $y$. Clearly, $x$ and $y$ are not adjacent in $\Gamma_i$ if and only if  $(x,y)\in S_k$, for some $k \neq i$, and  $p^k_{ii} \neq 0$, $k\neq i$, is equivalent to having an $S_i$-path of length two from $x$ to $y$, for any $(x,y)\in S_k$. Therefore, for a fixed $i$, only the values $k\neq i$ such that $ p^k_{ii}=0$ are to be considered.

			Suppose  $i=1$. Let $(x,y)\in S_2$. Then $\cV^{(x,y)}_{11}=\{z\}$ (since $p^2_{11}=1$). So there is a unique $S_1$-path $xzy$ from $x$ to $y$. Let $(x,y) \in S_3$. Then, $\cV^{(x,y)}_{21}=\{z\}$ and  $\cV^{(x,z)}_{11}=\{z'\}$ (since $p^3_{21}=1=p^2_{11}$); it follows that $xz' z y$ is the desired $S_1$-path.  Let $(x,y) \in S_4$. Since $p^4_{31}\neq 0$, from the previous arguments we may conclude that there is an $S_1$-path of length four from $x$ to $y$.

			The connectedness of the basis graph $\Gamma_i=(\Omega,S_i)$, for $i=2,3,4$, is proved by using very similar arguments.
		\end{proof}		
\begin{remark}
	{\em The basis graph $(\Omega, S_1)$ is the incidence graph of the dual of the double of $\cS$, that is the geometry $2\cS=(\cP\cup \cL,\Omega,\in)$; see \cite[p.2]{vm}. Since the fusion $(\Omega,\tilde\cR)$ exists only if $s=t$, then $(2\cS)^D$ is a weak generalized octagon of order $(s,1)$ \cite[p.21]{vm}; we refer the reader to \cite{vm} for additional information on weak generalized octagons.
}
		\end{remark}
		\subsection{Reconstructing the generalized quadrangle from the fusion} \label{sec_4}
		
		Let $\cY'=(\Omega',\tilde\cR')$ be an association scheme algebraically isomorphic to $\cY=(\Omega,\tilde\cR)$ via the isomorphism $\phi$ such that $S'_{i}=\phi(S_i)$, for $i=0,\ldots,4$.

		Our aim is to reconstruct a generalized quadrangle with parameters $(s,s)$ from $\cY'$. 
		
		From now on, ``clique'' will stand for  ``maximal $\{0, 1\}$-clique''. 
		
		\begin{lemma}
			For any $x\in\Omega'$, the set $S'_{1}(x)=\{y\in \Omega':(x,y)\in S'_{1}\}$ is partitioned in two cliques. 
		\end{lemma}
		\begin{proof}
			Let $y\in S'_1(x)$. Then there are  $p^1_{11}=s-1$ vertices 1-related to both $x$ and $y$. If $s=1$ (that is  $p^1_{11}=0$) or $s=2$ (that is  $p^1_{11}=1$), the result is clear. \\
			Let $s \geq 3$, and $u,v$ be two distinct vertices in $\cV^{(x,y)}_{11}$. Since $p^k_{11}\neq0$ only for $k=1,2$, we see that either $(u,v)\in S'_1$ or $(u,v)\in S'_2$. Assume the latter case occurs. Then $x,y\in \cV^{(u,v)}_{11}$. But this is a contradiction as $p^2_{11}=1$.
		\end{proof}
		
		Let $\cC$ be the set of all maximal $\{0,1\}$-cliques in $\cY'$. By counting in two ways the pairs $(x,C)$ with $x\in \Omega'$ and $C$ a clique on $x$, we see that $|\cC|=2(s+1)(s^2+1)$, which is precisely twice the number of points (of lines) of a GQ of order $(s,s)$.
		
		In light of the previous result, the idea is to select  $(s+1)(s^2+1)$ elements in $\cC$ (one for every $x\in\Omega'$) in such a way that these will be the  points of a (hypothetical) GQ. Clearly, the remaining cliques will be the lines.
		
		We will split the set $\cC$ in two disjoint subsets $\widehat \cP$ (points) and   $\widehat \cL$ (lines), each of size $(s+1)(s^2+1)$.

		\begin{lemma}\label{lem_1}
			Let  $C_1$ and $C_2$ be the two cliques on a vertex $x\in\Omega'$. Then, for every $u\in C_1\setminus\{x\}$ and $v\in C_2\setminus\{x\}$, $(u,v)\in S_{2}$ holds. 
		\end{lemma}
		\begin{proof}
			We have $v\in\cV^{(x,u)}_{1 i}=\{z\in \Omega':(x,z)\in S'_1,(z,u)\in S'_i\}$, for some  $i$ such that $p^{1}_{1i}$ is non-zero.  By looking at the matrix $M_{1}$, we see that $i\in\{0,1,2\}$. On the other hand, $\cV^{(x,u)}_{10}=\{u\}$, and $\cV^{(x,u)}_{1 1}=C_1\setminus\{x,u\}$, since $p^{1}_{11}=s-1$. Therefore, $v\in\cV^{(x,u)}_{1 2}$, that is $(u,v)\in S_{2}$, and $\cV^{(x,u)}_{1 2}=C_2\setminus\{x\}$, since $p^{1}_{12}=s$.
		\end{proof}
		\begin{lemma}\label{lem_3}
			Let $C=\{x_0,x_1,\ldots,x_{s}\}\in \cC$. For any $x_i\in C$, denote with $C'_i$ the clique on $x_i$ different from $C$. Then the  cliques $C'_i$, $i=0,1,\ldots,s$, are pairwise disjoint.
		\end{lemma}
		\begin{proof}  Let  $x_i,x_j$  distinct  vertices of $C$, and  $z\in C'_i\cap C'_j\neq \emptyset$. Then $C'_i$ and $C'_j$ are two cliques on $z$. Since $x_i,x_j\in C$, then $(x_i,x_j)\in S'_1$. By Lemma \ref{lem_1}, this yields that $x_i,x_j$ are in the same clique through  $z$, a contradiction. 
		\end{proof}
		
		Pick a vertex  $x_0\in \Omega'$. The idea is to split the vertices of $\Omega'$ into subsets, which we call {\em levels}, by considering  the distance between $x_0$ and the vertices of the given clique in the basis graph $(\Omega', S'_1)$. During this process, we will also ``label'' every clique by using the symbols $P$ and $L$. 
		
		\fbox{{\bf Level $\Lambda_0(x_0)$:}} We set $\Lambda_0(x_0)=\{x_0\}$;  it is obvious that $\Lambda_0(x_0)=S'_0(x_0)$.
		
		\fbox{{\bf Level $\Lambda_1(x_0)$:}} We denote the two cliques on $x_0$ by $P(x_0)$ and $L(x_0)$. 
		We use $\Lambda_1(x_0)$ to indicate the set of  the vertices of $P(x_0)\setminus\{x_0\}$ and  the vertices of $L(x_0)\setminus\{x_0\}$. It is obvious  that $|\Lambda_1(x_0)|=2s$ and  $\Lambda_1(x_0)=S'_1(x_0)$.
		
		\fbox{{\bf Level $\Lambda_2(x_0)$:}} For any vertex $x_1\in P(x_0)\setminus\{x_0\}$, denote the clique on $x_1$ different from $P(x_0)$ by $L(x_0,x_1)$. We set  $P(x_1)=P(x_0)$ and $L(x_1)=L(x_0,x_1)$.
		
		Set $\cL_2(x_0)=\{L(x_0,x_1): x_1\in P(x_0)\setminus\{x_0\}\}$. Clearly,  $|\cL_2(x_0)|=s$.

		\begin{corollary}\label{lem_2}
			For any $x_1\in P(x_0)\setminus\{x_0\}$, the vertices of $L(x_0,x_1)\setminus\{x_1\}$  are $2$-related to $x_0$. 
		\end{corollary}
		\begin{proof}
			We apply Lemma \ref{lem_1}.
		\end{proof}
		\begin{lemma}\label{lem_6}
			Let $x_2\in L(x_0,x_1)\setminus\{x_1\}$ and $x'_2\in L(x_0,x'_1)\setminus\{x'_1\}$, with $x_1,x'_1$ distinct vertices of $P(x_0)\setminus\{x_0\}$. Then $(x_2,x'_2)\in S'_3$.
		\end{lemma}
		\begin{proof}
			By Corollary \ref{lem_2}, $(x_0,x_2),(x_0,x'_2)\in S'_2$, so $x_2x_0x'_2$ is an $(S'_2,S'_2)$-path  from $x_2$ to $x'_2$. Therefore,  $x'_2\in \cV^{(x_0,x_2)}_{2i}$, for some $i$ such that $p^2_{2i}\neq 0$. By looking at the matrix $M_{2}$, we see that $i\in \{1,3,4\}$. On the other hand, $\cV^{(x_0,x_2)}_{21}=L(x_0,x_1)\setminus\{x_1,x_2\}$, since $p^{2}_{21}=s-1$. Assume $(x_2,x'_2)\in S'_4$. Then, by Lemma \ref{lem_1}, $x'_1\in \cV^{(x_2,x'_2)}_{21}$, with $p^4_{21}=0$; a contradiction. 
			It follows, $x'_2\in \cV^{(x_0,x_2)}_{23}$, and  $\cV^{(x_0,x_2)}_{23}=\bigcup\limits_{x'_1\in P(x_0)\setminus\{x_0,x_1\}}{L(x_0,x'_1})$, since $p^{2}_{23}=s(s-1)$.
		\end{proof}
		\comment{
			\begin{lemma}\label{lem_7}
				Let $x_2\in L(x_0,x_1)\setminus\{x_1\}$ and $y_1\in L(x_0)\setminus\{x_0\}$. Then $(x_2,y_1)\in S'_3$.
			\end{lemma}
			\begin{proof}
				By Corollary \ref{lem_2}, $(x_0,x_2)\in S'_2$, so $x_2x_0y_1$ is an $(S'_2,S'_1)$-path from $x_2$ to $y_1$. Therefore,  $y_1\in \cV^{(x_0,x_2)}_{1i}$, for some $i$ such that $p^2_{1i}\neq 0$. By looking at the matrix $M_{2}$ we see that $i\in \{1,2,3\}$. On the other hand, $\cV^{(x_0,x_2)}_{11}=\{x_1\}$ (since $p^{2}_{11}=1$) and $\cV^{(x_0,x_2)}_{12}=P(x_0)\setminus\{x_0,x_1\}$ (since $p^{2}_{12}=s-1$). Therefore, $y_1\in \cV^{(x_0,x_2)}_{13}$, and  $\cV^{(x_0,x_2)}_{13}=L(x_0)\setminus\{x_0\}$ (since $p^{2}_{13}=s$).
			\end{proof}
		}
		\begin{proposition}\label{prop_13}
			The cliques in $\cL_2(x_0)$ are pairwise disjoint. 
			Therefore, the vertices in $\cL_2(x_0)$ not in $P(x_0)$ are $s^2$.
		\end{proposition}
		\begin{proof}
			We apply Lemma \ref{lem_3}.
		\end{proof}
		For any vertex $y_1\in L(x_0)\setminus\{x_0\}$ denote with $P(x_0,y_1)$ the clique on $y_1$ different from $L(x_0)$. Then we set  $L(y_1)=L(x_0)$ and $P(y_1)=P(x_0,y_1)$. 
		
		Set $\cP_2(x_0)=\{P(x_0,y_1): y_1\in L(x_0)\setminus\{x_0\}\}$. Clearly,  $|\cP_2(x_0)|=s$. 
		
		\begin{remark}
		{\em In the hypothetical GQ, $P(x_0)$ will be a fixed point and  $L(x_0)$ a fixed line incident with $P(x_0)$. The cliques in $\cL_2(x_0)$ will be  the $s$ lines through $P(x_0)$ different from $L(x_0)$, while the cliques in $\cP_2(x_0)$ will be  the $s$ points on the line $L(x_0)$ different from $P(x_0)$.}
		\end{remark}

		By applying the same arguments as we did for $\cL_2(x_0)$, we can prove the following results.

		\begin{corollary}\label{lem_2bis}
			For any $y_1\in L(x_0)\setminus\{x_0\}$, the vertices of $P(x_0,y_1)\setminus\{y_1\}$  are $2$-related to $x_0$. 
		\end{corollary}
		\begin{lemma}\label{lem_6bis}
			Let $y_2\in P(x_0,y_1)\setminus\{y_1\}$ and $y'_2\in P(x_0,y'_1)\setminus\{y'_1\}$, with $y_1,y'_1$ distinct vertices of $L(x_0)\setminus\{x_0\}$. Then $(y_2,y'_2)\in S'_3$.
		\end{lemma}
		\comment{
			\begin{lemma}\label{lem_7bis}
				Let $y_2\in P(x_0,y_1)\setminus\{y_1\}$ and $x_1\in P(x_0)\setminus\{x_0\}$. Then $(y_2,x_1)\in S'_3$.
			\end{lemma}
		}
		\begin{proposition}\label{prop_5bis}
			The cliques in $\cP_2(x_0)$ are pairwise disjoint. 
			Therefore, the vertices in $\cP_2(x_0)$ not in $L(x_0)$ are $s^2$.
		\end{proposition}
		We refer to $\Lambda_2(x_0)$ as the set consisting of all vertices of the cliques in $\cP_2(x_0)\cup\cL_2(x_0)$ which are not in $\Lambda_1(x_0)$.
		\begin{remark}
			{\em Note that every clique in $\cL_2(x_0)$ is disjoint from $L(x_0)$; similarly, every clique in $\cP_2(x_0)$ is disjoint from $P(x_0)$.}
		\end{remark}
		\begin{proposition}\label{prop_5}
			The set $\Lambda_2(x_0)$ consists precisely of all the vertices which are $2$-related to $x_0$, that is $\Lambda_2(x_0)=S'_2(x_0)$. Hence, $|\Lambda_2(x_0)|=\eta_2=2s^2$.
		\end{proposition}
		\begin{proof}
			Take $L(x_0,x_1)\in \cL_2(x_0)$ and $P(x_0,y_1)\in \cP_2(x_0)$. Assume that $L(x_0,x_1)$ and $P(x_0,y_1)$ share a vertex $z$ different from $x_0$. By Lemma \ref{lem_1}, $(x_1,y_1)\in S'_2$. So $z \in \cV^{(x_1,y_1)}_{11}=\{x_0\}$, which implies $z=x_0$; a contradiction. This yields that $L(x_0,x_1)$ and $P(x_0,y_1)$ are disjoint. From Propositions \ref{prop_13} and \ref{prop_5bis},  we see that $|\Lambda_2(x_0)|=2s^2=\eta'_2$.
		\end{proof}
		\fbox{{\bf Level $\Lambda_3(x_0)$}:}
		For any $L(x_0,x_1)\in\cL_2(x_0)$ and  any $x_2\in L(x_0,x_1)\setminus\{x_1\}$, we denote  the clique on $x_2$ different from  $L(x_0,x_1)$ by $P(x_0,x_1,x_2)$. We set $L(x_2)=L(x_0,x_1)$ and $P(x_2)=P(x_0,x_1,x_2)$. Also, $P(x_0,x_1,x_2)$ coincides $P(x_1,x_2)$ if we choose $x_1$ instead of $x_0$.
		
		Let $\cP_3(x_0)=\{P(x_0,x_1,x_2): x_1\in P(x_0)\setminus\{x_0\}, x_2\in L(x_0,x_1)\setminus\{x_1\}\}$.  
		
		\begin{proposition}\label{prop_12}
			$|\cP_3(x_0)|=s^2$.  
		\end{proposition} 
		\begin{proof}
			This follows from Proposition \ref{prop_13}.  
		\end{proof}
		
		\begin{remark}
		{\em 	In the hypothetical GQ, the cliques in $\cP_3(x_0)$ will be the points collinear with the point $P(x_0)$ not incident with $L(x_0)$. For any fixed $x_1\in P(x_0)\setminus\{x_0\}$, the  cliques $P(x_0,x_1,x_2)$, $x_2\in L(x_0,x_1)$, will be  the $s$ points incident with the line $L(x_0,x_1)$ and different from $P(x_0)$.}
		\end{remark}
		\begin{lemma}\label{lem_4}
			For any $x_1\in P(x_0)\setminus\{x_0\}$ and $x_2\in L(x_0,x_1)\setminus\{x_1\}$, the vertices of $P(x_0,x_1,x_2)\setminus\{x_2\}$  are $3$-related to $x_0$. 
		\end{lemma}
		\begin{proof}
			Take $x_3\in P(x_0,x_1,x_2)\setminus\{x_2\}$. By Corollary \ref{lem_2}, $(x_2,x_0)\in S'_2$, so $x_3x_2x_0$ is an $(S_{1}, S_{2})$-path. Therefore,  $x_3\in\cV^{(x_2,x_0)}_{1i}$ for some  $i$ such that $p^{2}_{1i}$ is non-zero. 
			By looking at the matrix $M_{2}$, we see that $i\in \{1,2,3\}$. On the other hand, $\cV^{(x_2,x_0)}_{11}=\{x_1\}$ (since $p^{2}_{11}=1$), and  $\cV^{(x_2,x_0)}_{12}=L(x_0,x_1)\setminus\{x_1,x_2\}$ (since $p^{2}_{12}=s-1$). Therefore  $x_3\in \cV^{(x_2,x_0)}_{13}$, and  $\cV^{(x_2,x_0)}_{13}=P(x_0,x_1,x_2)\setminus\{x_2\}$ (since $p^{2}_{13}=s$). 
		\end{proof}
		\comment{
			\begin{lemma}\label{lem_12}
				Let $x_3\in P(x_0,x_1,x_2)\setminus\{x_2\}$ and $x'_1\in P(x_0)\setminus\{x_1\}$. Then  $(x_3,x'_1)\in S'_3$. {(\bf serve?)}
			\end{lemma}
			\begin{proof}
				By the  Lemma \ref{lem_1}, $(x_1,x_3)\in S'_2$, so $x_3x_1x'_1$ is an $(S'_2,S'_1)$-path from $x_3$ to $x'_1$. Therefore,  $x_3\in \cV^{(x_1,x'_1)}_{2i}$, for some $i$ such that $p^1_{2i}\neq 0$. By looking at the matrix $M_{1}$ we see that $i\in \{1,2,3\}$. On the other hand, $\cV^{(x_1,x'_1)}_{21}=L(x_0,x_1)\setminus\{x_1\}$ (since $p^{1}_{21}=s$) and $\cV^{(x_1,x'_1)}_{22}=\bigcup_{x\in P(x_0)\setminus\{x_1,x'_1\}}{(L(x)\setminus\{x\})}$ (since $p^{1}_{22}=s(s-1)$). Therefore, $x_3\in \cV^{(x_1,x'_1)}_{23}$, and  $\cV^{(x_1,x'_1)}_{23}=\bigcup_{x\in L(x_0,x_1)\setminus\{x_2\}}{(P(x)\setminus\{x\})}$ (since $p^{2}_{23}=s^2$).
			\end{proof}
		}
		\comment{
			\begin{lemma}
				Let $x_3\in P(x_0,x_1,x_2)$ and $y_1\in L(x_0)$. Then $(x_2,y_1)\in S'_4$.
			\end{lemma}
			\begin{proof}
				By Corollary \ref{lem_4}, $(x_0,x_3)\in S'_3$, so $x_3x_0y_1$ is an $(S'_3,S'_1)$-path from $x_3$ to $y_1$. Therefore,  $y_1\in \cV^{(x_0,x_3)}_{1i}$, for some $i$ such that $p^3_{1i}\neq 0$. By looking at the matrix $M_{3}$ we see that $i\in \{2,3,4\}$. On the other hand, $\cV^{(x_0,x_3)}_{12}=\{x_1\}$ (since $p^{3}_{12}=1$) and $\cV^{(x_0,x_3)}_{13}=P(x_0)\{x_0,x_1\}$ (since $p^{3}_{13}=s-1$). Therefore, $y_1\in \cV^{(x_0,x_3)}_{14}$, and  $\cV^{(x_0,x_4)}_{13}=L(x_0)\setminus\{x_0\}$ (since $p^{3}_{14}=s$).
			\end{proof}
		}

		\begin{corollary}\label{cor_8}
			Let $x_3\in P(x_0,x_1,x_2)\setminus\{x_2\}$ and  $x'_3\in P(x_0,x_1,x'_2)\setminus\{x'_2\}$, for $x_2,x'_2$ distinct vertices of $L(x_0,x_1)$, $x_1\in P(x_0)\setminus\{x_0\}$. Then $(x_3,x'_3)\in S_{3}$. 
		\end{corollary}
		
		\begin{proof}
			This follows by Lemma \ref{lem_6bis} applied to the cliques $P(x_1,x_2)\setminus\{x_2\}$ and $P(x_1,x'_2)\setminus\{x'_2\}$.
		\end{proof}
		%
		%
		%
		\begin{lemma}\label{lem_9}
			Let  $x_3\in P(x_0,x_1,x_2)$ and  $x'_3\in P(x_0,x'_1,x'_2)$, for $x_1,x'_1$ distinct vertices of $P(x_0)$. 
			Then $(x_3,x'_2),(x'_3,x_2)\in S'_4$ and $(x_3,x'_3)\in S'_3\cup S'_4$.
		\end{lemma}
		\begin{proof}
			By Lemma \ref{lem_6}, we have  $(x'_2,x_2)\in S'_3$, and  $x'_3\in\cV^{(x'_2,x_2)}_{1i}$ for some  $i$ such that $p^{3}_{1i}$ is non-zero. By looking at the matrix $M_{3}$,we see that $i\in\{2,3,4\}$. On the other hand, by Lemma \ref{lem_1},
			$\cV^{(x'_2,x_2)}_{12}=\{x'_1\}$, since $p^{3}_{12}=1$.  By applying Lemma \ref{lem_6}, we see that $\cV^{(x'_2,x_2)}_{13}=L(x_0,x'_1)\setminus\{x'_1,x'_2\}$, since $p^{3}_{13}=s-1$. Therefore, $(x'_3,x_2)\in S'_4$, and $\cV^{(x'_2,x_2)}_{14}=P(x_0,x'_1,x'_2)\setminus\{x'_2\}$, since $p^{3}_{14}=s$.
			
			Furthermore, $x_3\in \cV^{(x_2,x'_3)}_{1i}$ for some $i$ such that $p^4_{1i}\neq 0$. By looking at the matrix $M_4$, we see that $i\in\{3,4\}$.
		\end{proof}
		
		\begin{proposition}\label{prop_10}
			The cliques in $\cP_3(x_0)$ are pairwise disjoint. 
			Therefore, the vertices of the cliques in $\cP_3(x_0)$ which are not vertices of cliques in $\cL_2(x_0)$ are $s^3$.
		\end{proposition}  
		\begin{proof} 
			This follows from Lemma \ref{lem_9} and Proposition \ref{prop_12}.
		\end{proof} 

		\begin{remark}\label{rem_3}
			{\em By Lemmas \ref{lem_2} and \ref{lem_9}, any vertex of a  clique $P(x_0,x_1,x_2)$ in $\cP_3(x_0)$ which is not a vertex of a clique in $\cL_2(x_0)$ is in $\cV^{(x_0,x_1)}_{32}$. Since $|\cV^{(x_0,x_1)}_{32}|=p^1_{32}=s^2$, then $\cV^{(x_0,x_1)}_{32}$ consists precisely of all the vertices of the cliques $P(x_0,x_1,x_2)$ with $x_2\in L(x_0,x_1)\setminus\{x_1\}$.}
		\end{remark}

		For any $P(x_0,y_1)\in\cP_2(x_0)$ and  any $y_2\in P(x_0,y_1)\setminus\{y_1\}$, we denote the clique on $y_2$ different from  $P(x_0,y_1)$ by $L(x_0,y_1,y_2)$. We set $P(y_2)=P(x_0,y_1)$ and $L(y_2)=L(x_0,y_1,y_2)$. Also, $L(x_0,y_1,y_2)$ coincides with $L(y_1,y_2)$ if we choose $y_1$ instead of $x_0$.

		Let $\cL_3(x_0)=\{L(x_0,y_1,y_2): y_1\in L(x_0)\setminus\{x_0\}, y_2\in P(x_0,y_1)\setminus\{y_1\}\}$.  
		
		\begin{proposition}
			$|\cL_3(x_0)|=s^2$. 
		\end{proposition}
		\begin{proof}
			This follows from Proposition \ref{prop_5bis}.
		\end{proof} 

	\begin{remark}
		{\em	In the hypothetical GQ, the cliques in $\cL_3(x_0)$ will be the lines intersecting $L(x_0)$ not in $P(x_0)$. For any fixed $y_1\in L(x_0)\setminus\{x_0\}$, the  cliques $L(x_0,y_1,y_2)$, $y_2\in P(x_0,y_1)$, will be  the $s$ lines on the point $P(x_0,y_1)$ different from $L(x_0)$.}
	\end{remark}
		
		By applying the same arguments as we did for $\cP_3(x_0)$, we can prove the following results.

		\begin{lemma}\label{lem_4bis}
			For any $y_1\in L(x_0)\setminus\{x_0\}$ and $y_2\in P(x_0,y_1)\setminus\{y_1\}$, the vertices of $L(x_0,y_1,y_2)\setminus\{y_2\}$  are $3$-related to $x_0$. 
		\end{lemma}
		\comment{
			\begin{lemma}
				Let $y_3\in L(x_0,y_1,y_2)\setminus\{y_2\}$ and $y'_1\in L(x_0)\setminus\{y_1\}$. Then  $(y_3,y'_1)\in S'_3$. {(\bf serve?)}
			\end{lemma}
		}
		\comment{
			\begin{lemma}\
				Let $x_3\in P(x_0,x_1,x_2)$ and $y_1\in L(x_0)$. Then $(x_2,y_1)\in S'_4$.
			\end{lemma}
			\begin{proof}
				By Corollary \ref{lem_4}, $(x_0,x_3)\in S'_3$, so $x_3x_0y_1$ is an $(S'_3,S'_1)$-path from $x_3$ to $y_1$. Therefore,  $y_1\in \cV^{(x_0,x_3)}_{1i}$, for some $i$ such that $p^3_{1i}\neq 0$. By looking at the matrix $M_{3}$ we see that $i\in \{2,3,4\}$. On the other hand, $\cV^{(x_0,x_3)}_{12}=\{x_1\}$ (since $p^{3}_{12}=1$) and $\cV^{(x_0,x_3)}_{13}=P(x_0)\{x_0,x_1\}$ (since $p^{3}_{13}=s-1$). Therefore, $y_1\in \cV^{(x_0,x_3)}_{14}$, and  $\cV^{(x_0,x_4)}_{13}=L(x_0)\setminus\{x_0\}$ (since $p^{3}_{14}=s$).
			\end{proof}
		}

		\begin{corollary}\label{cor_8bis}
			Let $y_3\in L(x_0,y_1,y_2)\setminus\{y_2\}$ and  $y'_3\in L(x_0,y_1,y'_2)\setminus\{y'_2\}$, for $y_2,y'_2$ distinct vertices of $P(x_0,y_1)$, $y_1\in L(x_0)\setminus\{x_0\}$. Then $(y_3,y'_3)\in S_{3}$. 
		\end{corollary}

		%
		%
		%
		\begin{lemma}\label{lem_9bis}
			Let  $y_3\in L(x_0,y_1,y_2)$ and  $y'_3\in L(x_0,y'_1,y'_2)$, for $y_1,y'_1$ distinct vertices of $L(x_0)$. 
			Then $(y_3,y'_2),(y'_3,y_2)\in S'_4$ and $(y_3,y'_3)\in S'_3\cup S'_4$.
		\end{lemma}
		\begin{proposition}\label{prop_10bis}
			The cliques in $\cL_3(x_0)$ are pairwise disjoint. 
			Therefore, the vertices of the cliques in $\cL_3(x_0)$ which are not vertices of cliques in $\cP_2(x_0)$ are $s^3$.
		\end{proposition}  
		\begin{proof} 
			This follows from Lemma  \ref{lem_9bis} and Proposition \ref{prop_5bis}.
		\end{proof} 
		We refer to $\Lambda_3(x_0)$ as the set consisting of all vertices of the cliques in $\cP_3(x_0)\cup\cL_3(x_0)$ which are not in $\Lambda_2(x_0)$. 
		\begin{proposition}\label{prop_5tris}
			The set $\Lambda_3(x_0)$ consists precisely of all the vertices $3$-related to $x_0$, that is $\Lambda_3(x_0)=S'_3(x_0)$. Hence, $|\Lambda_3(x_0)|=\eta_3 =2s^3$.
		\end{proposition}
		\begin{proof}
			Take $x_3\in P(x_0,x_1,x_2)\in \cP_3(x_0)$ and $y_3\in L(x_0,y_1,y_2)\in \cL_3(x_0)$. We now show that $(y_1,x_3)\in S'_4$.  By Lemma \ref{lem_4}, applied to $y_1\in L(x_0)=L(y_1)$ and $x_2\in L(x_0,x_1)=L(y_1,x_0,x_1)$, we have $(x_2,y_1)\in S'_3$. Hence, $x_3x_2y_1$ is an $(S'_1,S'_3)-$path from $x_3$ to $y_1$. Therefore,  $x_3\in\cV^{(x_2,y_1)}_{1i}$, for some  $i$ such that $p^{3}_{1i}$ is non-zero. 
			
			By looking at the matrix $M_{3}$, we see that $i\in \{2,3,4\}$. On the other hand, by Lemma \ref{lem_1}, $\cV^{(x_2,y_1)}_{12}=\{x_1\}$, since $p^{3}_{12}=1$, and  $\cV^{(x_2,y_1)}_{13}=L(x_0,x_1)\setminus\{x_1,x_2\}$, by Lemma \ref{lem_4bis} applied to $L(y_1,x_0,x_1)$, since $p^{3}_{13}=s-1$. Therefore  $x_3\in \cV^{(x_2,y_1)}_{14}$, and  $\cV^{(x_2,y_1)}_{14}=P(x_0,x_1,x_2)\setminus\{x_2\}$, since $p^{3}_{14}=s$. 
			
			By Lemma \ref{lem_2}, $x_3y_1y_3$ is an $(S'_4,S'_2)$-path from $x_3$ to $y_3$. Hence,  $x_3\in\cV^{(y_1,y_3)}_{4i}$, for some  $i$ such that $p^{2}_{4i}$ is non-zero. By looking at the matrix $M_{2}$, we see that $i\in \{2,3,4\}$. This implies that $(x_3,y_3)\not\in S'_0$.  So  $P(x_0,x_1,x_2)$ and $ L(x_0,y_1,y_2)$ are disjoint.
			
			From  Propositions \ref{prop_10} and \ref{prop_10bis}, we see that $|\Lambda_3(x_0)|=2s^3=\eta'_3$.
		\end{proof}
		
		\begin{remark}\label{rem_1bis}
			{\em By Lemmas \ref{lem_2} and \ref{lem_4bis}, any vertex of a  clique $L(x_0,y_1,y_2)$ in $\cL_3(x_0)$ which is not a vertex of a clique in $\cL_2(x_0)$ is in $\cV^{(x_0,y_1)}_{32}$. Since $|\cV^{(x_0,y_1)}_{32}|=p^1_{32}=s^2$, $\cV^{(x_0,y_1)}_{32}$ consists precisely of all the vertices of the cliques $L(x_0,y_1,y_2)$ with $y_2\in P(x_0,y_1)\setminus\{y_1\}$.}
		\end{remark}

		\fbox{{\bf Level $\Lambda_4(x_0)$}:}
		For any $P(x_0,x_1,x_2)\in\cP_3(x_0)$ and  any $x_3\in P(x_0,x_1,x_2)\setminus\{x_2\}$, we denote by $L(x_0,x_1,x_2,x_3)$ the clique on $x_3$ different from  $P(x_0,x_1,x_2)$. We set $P(x_3)=P(x_0,x_1,x_2)$ and $L(x_3)=L(x_0,x_1,x_2,x_3)$. Furthermore,  $L(x_0,x_1,x_2,x_3)=L(x_1,x_2,x_3)$ if we choose $x_1$ instead of $x_0$, and $L(x_0,x_1,x_2,x_3)=L(x_2,x_3)$ if we choose $x_2$ instead of $x_0$.
		
		Let $\cL_4(x_0)=\{L(x_0,x_1,x_2,x_3): x_1\in P(x_0)\setminus\{x_0\}, x_2\in L(x_0,x_1)\setminus\{x_1\},x_3\in P(x_0,x_1,x_2)\setminus\{x_2\}\}$. 
		\begin{proposition}
			$|\cL_4(x_0)|=s^3$.   
		\end{proposition}
		\begin{proof}
			This follows from Proposition \ref{prop_10}.
		\end{proof}
		
		\begin{remark}
		{\em	In the hypothetical GQ, the cliques in $\cL_4(x_0)$ will be the lines not incident with $P(x_0)$ and intersecting some line through $P(x_0)$.}
		\end{remark}
		
		\begin{lemma}\label{lem_10}
			For any $x_1\in P(x_0)\setminus\{x_0\}$, $x_2\in L(x_0,x_1)\setminus\{x_1\}$ and $x_3\in P(x_0,x_1,x_2)\setminus\{x_1,x_2\}$, the vertices of $L(x_0,x_1,x_2,x_3)\setminus\{x_3\}$  are $4$-related to $x_0$.
		\end{lemma}
		\begin{proof}
			Take $x_4\in L(x_0,x_1,x_2,x_3)\setminus\{x_3\}$. By Lemma \ref{lem_4}, $(x_0,x_3)\in S'_3$, so $x_4x_3x_0$ is an $(S_{1}, S_{3})$-path. Therefore,  $x_4\in\cV^{(x_3,x_0)}_{1i}$ for some  $i$ such that $p^{3}_{1i}$ is non-zero. 
			By looking at the matrix $M_{3}$, we see that $i\in \{2,3,4\}$. On the other hand, $\cV^{(x_3,x_0)}_{12}=\{x_2\}$ by Corollary \ref{lem_2} (since $p^{3}_{12}=1$), and  $\cV^{(x_3,x_0)}_{13}=P(x_0,x_1,x_2)\setminus\{x_2,x_3\}$ by Lemma \ref{lem_4} (since $p^{3}_{13}=s-1$). Therefore  $x_4\in \cV^{(x_3,x_0)}_{14}$, and  $\cV^{(x_3,x_0)}_{14}=L(x_0,x_1,x_2,x_3)\setminus\{x_3\}$ (since $p^{3}_{14}=s$). 
		\end{proof}
		%
		%
		
		%
		\begin{corollary}\label{cor_9}
			Let $x_4\in L(x_0,x_1,x_2,x_3)$ and  $x'_4\in L(x_0,x_1,x_2,x'_3)$, with $x_3,x'_3$ distinct vertices of $P(x_0,x_1,x_2)$, $x_1\in P(x_0)\setminus\{x_0\}$, $x_2\in L(x_0,x_1)\setminus\{x_1\}$. Then $(x_4,x'_4)\in S'_3$.
		\end{corollary}
		\begin{proof}
			This follows from Lemma \ref{lem_6} applied to the cliques $L(x_2,x_3)\setminus\{x_3\}$ and $L(x_2,x'_3)\setminus\{x'_3\}$.
		\end{proof}
		\begin{corollary}
			For any fixed $x_1\in P(x_0)\setminus\{x_0\}$, $x_2\in L(x_0,x_1)\setminus\{x_1\}$, the cliques $L(x_0,x_1,x_2,x_3)$ and  $L(x_0,x_1,x_2,x'_3)$, with $x_3,x'_3$ distinct vertices of $P(x_0,x_1,x_2)\setminus\{x_2\}$,  are pairwise disjoint.  
		\end{corollary}
		\begin{proof}
			It immediately follows from Lemma \ref{lem_3}.
		\end{proof}
		\begin{corollary}\label{cor_29}
			Let $x_4\in L(x_0,x_1,x_2,x_3)$ and  $x'_4\in L(x_0,x_1,x'_2,x'_3)$, for $x_2,x'_2$ distinct vertices of $L(x_0,x_1)$, $x_1\in P(x_0)\setminus\{x_0\}$. Then $(x_4,x'_3),(x'_4,x_3)\in S'_4$ and $(x_4,x'_4)\in S'_3\cup S'_4$. 
		\end{corollary}
		\begin{proof}
			This follows from Lemma \ref{lem_9bis} applied to the cliques $L(x_1,x_2,x_3)\setminus\{x_3\}$ and $L(x_1,x'_2,x'_3)\setminus\{x'_3\}$.
		\end{proof}

		\begin{corollary}
			For any fixed $x_1\in P(x_0)\setminus\{x_0\}$ and $x_2,x'_2$ distinct vertices of $L(x_0,x_1)\setminus\{x_1\}$, the cliques $L(x_0,x_1,x_2,x_3)$ and  $L(x_0,x_1,x'_2,x'_3)$, with $x_3\in P(x_0,x_1,x_2)\setminus\{x_2\}$ and $x'_3\in P(x_0,x_1,x'_2)\setminus\{x'_2\}$,  are pairwise disjoint.  
		\end{corollary}  
		\begin{proof}
			It immediately follows from Corollary \ref{cor_29}.
		\end{proof}
		\begin{lemma}\label{lem_11}
			Let  $x_4\in L(x_0,x_1,x_2,x_3)$ and  $x'_4\in L(x_0,x'_1,x'_2,x'_3)$, for $x_1,x'_1$ distinct vertices of $P(x_0)$. 
			Then $(x_4,x'_4)\notin S'_0$.
		\end{lemma}
		\begin{proof}
			By Lemma \ref{lem_9}, we have  $(x'_3,x_3)\in S'_3\cup S'_4$, and  $x'_4\in\cV^{(x'_3,x_3)}_{1i}$ for some  $i$ such that $p^{k}_{1i}$ is non-zero, for $k\in\{3,4\}$. Assume $(x'_3,x_3)\in S'_3\cup S'_4$.  By looking at the matrices $M_{3}$ and $M_4$, we see that $i\in\{2,3,4\}$. Therefore, $x_4x_3x'_4$ is an $(S'_1,S'_i)$-path, for some $i\in\{2,3,4\}$. Hence, $x_4\in\cV^{(x_3,x'_4)}_{1j}$, with $(x_3,x'_4)\in S'_i$, for some $i\in\{2,3,4\}$, such that $p^i_{1j}\neq0$. By considering the second row of the matrices $M_i$, $i=2,3,4$, we see that $j\neq 0$. This proves the result.
		\end{proof}
		\begin{proposition}\label{prop_11}
			The cliques in $\cL_4(x_0)$ are pairwise disjoint. Therefore, the vertices of the cliques in $\cL_4(x_0)$ which are not vertices of cliques in $\cP_3(x_0)$ are $s^4=\eta'_4$. 
		\end{proposition}  
		\begin{proof} 
			This follows from Lemma \ref{lem_11} and Proposition \ref{prop_10}.
		\end{proof} 
		Let $\cP_4(x_0)=\{P(x_0,y_1,y_2,y_3): y_1\in L(x_0)\setminus\{x_0\}, y_2\in P(x_0,y_1)\setminus\{y_1\},y_3\in L(x_0,y_1,y_2)\setminus\{y_2\}\}$. We set $L(y_3)=L(x_0,y_1,y_2)$ and $P(y_3)=P(x_0,y_1,y_2,y_3)$. Furthermore,  $P(x_0,y_1,y_2,y_3)=P(y_1,y_2,y_3)$ if we choose $y_1$ instead of $x_0$, and $P(x_0,y_1,y_2,y_3)=P(y_2,y_3)$ if we choose $y_2$ instead of $x_0$.
		\begin{proposition}
			$|\cP_4(x_0)|=s^3$.   
		\end{proposition}
		
		\begin{proof}
			This follows from Proposition \ref{prop_10bis}.
		\end{proof}
		
		\begin{remark}
			{\em In the hypothetical GQ, the cliques in $\cP_4(x_0)$ will be the points not collinear with $P(x_0)$.}
		\end{remark}

		By applying the same arguments as we did for $\cL_4(x_0)$, we can prove the following results.
		
		\begin{lemma}\label{lem_10bis}
			For any $y_1\in L(x_0)\setminus\{x_0\}$, $y_2\in P(x_0,y_1)\setminus\{y_1\}$ and $y_3\in L(x_0,y_1,y_2)\setminus\{y_1,y_2\}$, the vertices of $P(x_0,y_1,y_2,y_3)\setminus\{y_3\}$  are $4$-related to $x_0$. 
		\end{lemma}

		\begin{corollary}\label{cor_9bis}
			Let $y_4\in P(x_0,y_1,y_2,y_3)$ and  $y'_4\in L(x_0,y_1,y_2,y'_3)$, with $y_3,y'_3$ distinct vertices of $L(x_0,y_1,y_2)$, $y_1\in P(x_0)\setminus\{x_0\}$, $y_2\in P(x_0,y_1)\setminus\{y_1\}$. Then $(y_4,y'_4)\in S'_3$. 
		\end{corollary}
		\begin{corollary}
			For any fixed $y_1\in L(x_0)\setminus\{x_0\}$, $y_2\in P(x_0,y_1)\setminus\{y_1\}$, the cliques $P(x_0,y_1,y_2,y_3)$ and  $P(x_0,y_1,y_2,y'_3)$, with $y_3,y'_3$ distinct vertices of $L(x_0,y_1,y_2)\setminus\{y_2\}$,  are pairwise disjoint.  
		\end{corollary}  
		\begin{corollary}\label{cor_29bis}
			Let $y_4\in P(x_0,y_1,y_2,y_3)$ and  $y'_4\in P(x_0,y_1,y'_2,y'_3)$, for $y_2,y'_2$ distinct vertices of $P(x_0,y_1)$, $y_1\in L(x_0)\setminus\{x_0\}$. Then $(y_4,y'_3),(y'_4,y_3)\in S'_4$ and $(y_4,y'_4)\in S'_3\cup S'_{4}$. 
		\end{corollary}
		\begin{corollary}
			For any fixed $y_1\in L(x_0)\setminus\{x_0\}$ and $y_2,y'_2$ distinct vertices of $P(x_0,y_1)\setminus\{y_1\}$, the cliques $P(x_0,y_1,y_2,y_3)$ and  $P(x_0,y_1,y'_2,y'_3)$, with $y_3\in L(x_0,y_1,y_2)\setminus\{y_2\}$ and $y'_3\in L(x_0,y_1,y'_2)\setminus\{y'_2\}$,  are pairwise disjoint.  
		\end{corollary}  
		\begin{lemma}\label{lem_11bis}
			Let  $y_4\in P(x_0,y_1,y_2,y_3)$ and  $y'_4\in P(x_0,y'_1,y'_2,y'_3)$, for $y_1,y'_1$ distinct vertices of $L(x_0)$. 
			Then $(y_4,y'_4)\in S'_2\cup S'_3\cup S'_4$.
		\end{lemma}
		\begin{proposition}\label{prop_11bis}
			The cliques in $\cP_4(x_0)$ are pairwise disjoint. Therefore, the vertices of the cliques in $\cP_4(x_0)$ which are not vertices of cliques in $\cL_3(x_0)$ are $s^4=\eta_4$.
		\end{proposition}  
		We refer to $\Lambda_4(x_0)$ as the set} consisting of all vertices of the cliques in $\cL_4(x_0)$ which are not vertices of $\cP_3(x_0)$.
	
	\begin{proposition}\label{prop_14}
		The set $\Lambda_4(x_0)$ coincides with the set of all vertices of the cliques in $\cP_4(x_0)$  which are not vertices of $\cL_3(x_0)$. Therefore, $|\Lambda_4(x_0)|=\eta_4=s^4$ and  $\Lambda_4(x_0)=S'_4(x_0)$. 
	\end{proposition}
	
	\begin{proof}
		This is an immediate consequence of Propositions \ref{prop_11} and \ref{prop_11bis}.
	\end{proof}
	Set 
	\[
	\widehat \cP=\cP_1(x_0)\cup\cP_2(x_0)\cup\cP_3(x_0)\cup\cP_4(x_0)
	\]
	and 
	\[
	\widehat \cL=\cL_1(x_0)\cup\cL_2(x_0)\cup\cL_3(x_0)\cup\cL_4(x_0).
	\]
	Since $\{S'_0(x_0),S'_1(x_0),S'_2(x_0),S'_3(x_0),S'_4(x_0)\}$ is a partition of the vertex set $\Omega'$ , we see that through every  $x\in\Omega'$ there is one clique in $\widehat\cP$, denoted by $P(x)$, and one clique in $\widehat\cL$, denoted by $L(x)$. Note also that the same sets of cliques $\widehat\cP$ and $\widehat\cL$ are constructed as we did before by using any vertex $x\in \Omega'$ instead of $x_0$. 
	
	We call the elements of $\widehat\cP$ \emph{points} and those of $\widehat\cL$ \emph{lines}.  We say that a point $P \in \widehat\cP$ and  a line $L \in \widehat\cL$ are \emph{incident}, and we will write $P\,{\rm \widehat I}\, L$, if $P$ and $L$ have a vertex in common.  
	
	We are going to show that  $\widehat{\cS}=(\widehat{\cP},\widehat{\cL},\  {\rm \widehat I})$ is a generalize quadrangle of order $(s,s)$. 
	
	Since every point has $s+1$ vertices, each of which is on a unique line, it follows that every point is incident with $s+1$ lines. Similarly, we find that every line is incident with $s+1$ points. So, from the maximality of the cliques, axioms (GQ1) and (GQ2) are satisfied. 
	\begin{theorem}
		Let $P\in\widehat\cP$ and $L\in\widehat\cL$ be not incident, i.e., $P$ and $L$ have no vertex in common.  The there exists a unique clique $Q\in\widehat\cP$ and a unique clique $M\in\widehat\cL$ such that $P\,{\rm \widehat I}\,M\,{\rm \widehat I}\,Q\,{\rm \widehat I}\,L$.
	\end{theorem}
	\begin{proof}
		Let $x\in P$ and $y\in L$. Different cases are treated separately depending on the relation where $(x,y)$ lies.
		
		Clearly $(x,y)\notin S'_0\cup S'_1$.
		
		Assume $(x,y)\in S'_2$. Since $p^2_{11}=1$, there exists a unique $z$ which is 1-related to $x$ and $y$. We set $M=L(z)=L(x)$ and $Q=P(z)=P(y)$. 
		
		Assume $(x,y)\in S'_3$. Since $p^3_{21}=1$, we have $\cV^{(x,y)}_{21}=\{z\}$. Note that $P(z)=P(y)$. Since $(x,z)\in S'_2$, we have either $P(z)\cap L(x)\neq \emptyset$ or $L(z)\cap P(x)\neq \emptyset$. Assume $P(z)\cap L(x)=\{v\}$.  Since $P(z)=P(y)$, then $v\in\cV^{(x,y)}_{11}$. On the other hand, $p^3_{11}=0$; a contradiction. Therefore, $L(z)\cap P(x)=\{v\}$, with $L(v)=L(z)$. In this case, we set $M=L(v)$ and $Q=P(y)$ to get the result.
		
		Assume $(x,y)\in S'_4$. Since $p^4_{31}=2$, we have $\cV^{(x,y)}_{31}=\{z,z'\}$. Suppose $z,z'\in L(y)$. Since $(x,z)\in S'_3$ and $p^3_{21}=1$, we have $\cV^{(x,z)}_{21}=\{v\}$. Note that $L(z)=L(z')=L(y)$ and $P(v)=P(z)$. This implies that there is precisely  one clique on $v$ and one clique on $x$ sharing a vertex. But $P(v)\cap L(x)\neq \emptyset$ is not possible as $p^3_{11}=0$. So, necessarily $L(v)\cap P(x)=\{w\}$. Note that we may write $P(x)=P(y,z,v,w)$. By applying the same arguments to $z'$, we obtain  vertices  $w'\neq w$  and $v'\neq v$ such that $P(x)=P(y,z',v',w')$. But this is not possible by Proposition \ref{prop_11}.
		By symmetry,  $z,z'\in P(y)$ cannot hold. Therefore, without loss of generality, we may assume $z\in L(y)$ and $z'\in P(y)$. 
		
		From the above arguments, we see that $M=L(v)$ and $Q=P(v)=P(z)$ are such that $P\,{\rm \widehat I}\,M\,{\rm \widehat I}\,Q\,{\rm \widehat I}\,L$.
		
		Since $z'\in P(y)$ with $(x,z')\in S'_3$, then $\cV^{(x,z')}_{21}=\{v'\}$. But $v'\in P(z')=P(y)$ cannot occur as $p^4_{21}=0$. So $v'\in L(z')$,  necessarily. Since $(x,v')\in S'_2$ and  $L(v')=L(z')$, we have either $P(v')\cap L(x)\neq\emptyset$ or $L(z')\cap P(x)\neq\emptyset$.  Since $(x,z')\in S'_3$ the latter case cannot occur. From this we get the uniqueness of $Q$ and $M$ such that $P\,{\rm \widehat I}\,M\,{\rm \widehat I}\,Q\,{\rm \widehat I}\,L$.
	\end{proof}

\end{document}